\documentclass[doublespacing]{elsart}

\usepackage{amssymb}
\usepackage{amsmath}
\newcommand{\F}{\ensuremath{\textup{F}}}
\newcommand{\GO}{\ensuremath{{\rm O}}}
\newcommand{\SO}{\ensuremath{{\rm SO}}}

\newcommand{\Spp}{\ensuremath{\textup{Sp}}}

\newcommand{\M}{\ensuremath{{\textup M}}}
\newcommand{\ZZ}{\ensuremath{{\textup Z}}}
\newcommand{\B}{\ensuremath{{\textup B}}}

\newcommand{\HH}{\ensuremath{{\textup H}}}

\newcommand{\T}{\ensuremath{{\textup T}}}

\newcommand{\V}{\ensuremath{{\textup V}}}
\newcommand{\U}{\ensuremath{{\textup U}}}

\newcommand{\Z}{\ensuremath{\mathbb{Z}}}
\newcommand{\A}{\ensuremath{\mathbb{A}}}
\newcommand{\R}{\ensuremath{\mathbb{R}}}
\newcommand{\C}{\ensuremath{\mathbb{C}}}
\newcommand{\Q}{\ensuremath{\mathbb{Q}}}
\newcommand{\Sc}{\ensuremath{\mathcal{S}}}
\newcommand{\OO}{\ensuremath{\mathcal{O}}}

\newcommand{\SL}{\ensuremath{\textup{SL}}}
\newcommand{\SLT}{\ensuremath{\widetilde{\textup{SL}}_{2}}}
\newcommand{\GL}{\ensuremath{\textup{GL}}}

\newcommand{\om}[2][]{\ensuremath{\Omega_{#1  \, #2 }}}

\newcommand{\rd}[1]{\ensuremath{\textup{ord}( #1 )}}
\renewcommand{\theequation}{\thesection.\arabic{equation}}

\begin{document}

\begin{frontmatter}


 \title{THETA CORRESPONDENCE \\OF AUTOMORPHIC CHARACTERS}
\author{Kobi Snitz}
\ead{snitz@math.bgu.ac.il}
\address{snitz@math.bgu.ac.il \\
 The Centre for Advanced Studies in Mathematics\\
  Ben Gurion University of the Negev.\\
P O Box 653 Be'er Sheva \\
84105 Israel.\\
phone number : 972-3-687-9485 \\
fax number: 972-8-6477648}




\pagebreak
\begin{abstract}
This paper describes the lifting of automorphic characters
of $\GO(3)(\A)$ to $\SLT(\A)$.  It does so by matching the image of
this
lift with the lift of automorphic characters from $\GO(1)(\A)$ to
$\SLT(\A)$. Our matching actually gives a matching of individual
automorphic forms, and not just of representation spaces.
Let $\V$ be a $3-$ dimensional quadratic vector space
and $\U$ a certain $1-$ dimensional quadratic space. To an automorphic
form $I_{\V}(\chi,\varphi)$ determined by the Schwartz function
$\varphi\in \Sc(\V(\A))$ in the lift of the character $\chi$ we
match an automorphic form $I_{\U}(\mu,\varphi_{0})$ determined by the
Schwartz function $\varphi_{0}\in \Sc(\U(\A))$ in the lift of the
character $\mu$.

Our work shows that,
the space $\U$ is explicitly determined by the character $\chi$. The
character $\mu$ is explicitly determined by the space $\V$ and the
function $\varphi_{0}$ is given by an orbital integral involving
$\varphi$.
\end{abstract}

\begin{keyword}
theta lift \sep automorohic characters
\PACS 
\end{keyword}
\end{frontmatter}
\pagebreak
\tableofcontents
\section{Introduction}

The Siegel-Weil formula 
identifies the integral of a certain theta function as 
an Eisenstein series on a symplectic group. In modern language, this is 
an example of a global theta correspondence for the dual pair 
${\rm Sp}(W),\GO(\V)$ in which the trivial representation of
$\GO(\V)(\A)$ corresponds to an automorphic representation of
$\Spp(W)(\A)$ generated by the values of Eisenstein series.
The group $\GO(\V)(\A)$ has many other one dimensional 
automorphic representations (automorphic characters) and it is
an open problem to describe where they all go in the correspondence.

This paper deals with the dual pair, $\SLT, \GO(3)$. The local and global 
versions of the theta correspondence, in this case, 
were studied by Waldspurger \cite{W1} \cite{W2} and he obtained  very
explicit information about corresponding representations. However
he did not consider global one dimensional representations of $\GO(3)(\A)$.

Our main results describes explicitly the lifts of automorphic
characters $\chi \neq 1$ of $\GO(3)$ in terms of theta series coming
from various $\GO(1)$'s, these results provide a complement to the
work of waldspurger and are analagues to the Seigel-Weil Formula. 

\subsection{Description of the result}
Let $\F$ be a number field and let $\B$ be a quaternion
algebra over $\F$. 
For any place $v$ of $\F$, let
$\F_{v}$ be the completion of $\F$,
and let $\B_{v}=\B\otimes\F_{v}$. The local algebra $\B_{v}$ can be
either split, $\B_{v}\simeq \M_{2}(\F_{v})$
or non-split, so that $\B_{v}$ is a division algebra. The 
invariant $\textup{Inv}(\B_{v})$ is
defined as $1$ when  $\B_{v}$ is split and $-1$ when $\B_{v}$ is 
non-split. 
$
\prod_{v}\textup{Inv}(\B_{v})=1.
$
Let 
\[
S=S_{\B}:=\{v: \textup{Inv}(\B_{v})=-1\}.
\]
By the product formula 
$
\prod_{v}\textup{Inv}(\B_{v})=1.
$
the cardinality of $S$, is even. The isomorphism class of
$\B$ is determined uniquely by $S$.

Let $\V$ be the trace zero elements of $\B$,
\[
\V := \{x \in \B : {\rm Tr}(x)=0\},
\]
and equip $\V$ with a
quadratic form $Q(x)=\nu(x)$, where $\nu:\B \rightarrow \F$ is the
reduced norm. The space $\V$ is thus a three dimensional
quadratic space over $\F$ and, up to scaling $Q$, all three dimensional
spaces over $\F$ are constructed this way. Since the dimension of $\V$ is  
odd, the orthogonal group factors as
$\GO(\V)=\SO(\V)\times \{\pm 1_{\V}\}$. 
Moreover, the group $\HH=\B^{\times}$ acts on $\V$ by conjugation and
the action gives rise to the exact sequence
\[
1\rightarrow Z\rightarrow \HH \rightarrow \SO(\V)\rightarrow 1,
\]
where $Z$ is the center of $\HH$. We can thus identify $\HH$ with 
${\rm GSpin}(\V)$

Let $\A=\A_{\F}$ be the adeles of $\F$ and $\A^{\times}$ be the ideles
of $\F$. Let $\chi$ be a character of $\A^{\times}/\F^{\times}$ such
that $\chi^{2}=1$. Such a character is 
given by 
\[
\chi(a)=(a,-\kappa)_{\A}=\prod_{v}(a_{v},-\kappa)_{\F_{v}}=
\prod_{v}\chi_{v}(a_{v})
\]
for some  $\kappa \in \F^{\times}$
where the coset $\kappa (\F^{\times})^{2}$ in
$\F^{\times}/(\F^{\times})^{2}$ is uniquely determined and where
the global Hilbert symbol $(\;,\;)_{\A}$ is defined by the local
Hilbert symbols $(\;,\;)_{\F_{v}}$.

Since $\nu:\F^{\times}\simeq Z \rightarrow \F^{\times}$ is the map
$x \mapsto x^{2}$ the composition $\chi\circ\nu$ is a quadratic character on
$\HH(\A)$ which is trivial on 
$Z(\A)\HH(\F)$. Hence it is an
automorphic character of $\HH(\A)$ which descends to $\SO(\A)$. All
automorphic characters of $\SO(\A)$ are obtained in this
way.

Let $\SLT(\A)$, be the $2$-fold metaplectic cover of $\SL(\A)$
Fix a nontrivial additive character $\psi:\A/\F \rightarrow \C^{1}$.
Then 
$\SLT(\A)$ acts on the Schwartz space 
$\Sc(\V(\A))$ by the Weil representation $\omega=\omega_{\psi}$
determined by $\psi$.
The action of $\SLT(\A)$ 
commutes with the linear action of $\HH(\A)$ on 
$\Sc(\V(\A))$ given by $h\cdot\varphi(x)=\varphi(h^{-1}\cdot x)$. Using
the commuting actions, we define for
$h\in \HH(\A), g\in \SLT(\A)$ and $\varphi \in \Sc(\V(\A))$
the theta function
\begin{equation}\label{tet}
\theta(g,h,\varphi) =\sum_{x \in  V(\F)}\omega(g)(h\cdot
\varphi)(x)= \sum_{x \in V(\F)}\omega(g)\varphi(h^{-1}\cdot x).
\end{equation}
This function
is left invariant by $\SL(\F)\times\HH(\F)$. 

We now assume that
$\B(\A)$ is a division algebra, i.e. $S\neq \emptyset$. In this case,
the orthogonal group $\SO(\V)$ is anisotropic and the quotient 
$\HH(\F)\ZZ(\A)\backslash \HH(\A)$ is compact. Thus, for a quadratic
character
$\chi$ as above, the integral
\begin{equation}\label{cor3}
I_{\V}(\chi,\varphi)(g) = \int\limits_{\HH(\F)\ZZ(\A)\backslash \HH(\A)} 
\theta(g,h,\varphi) \chi(\nu(h))dh
\end{equation}
is absolutely convergent and defines an automorphic form on
$\SLT(\A)$.

The space 
\begin{equation}\label{setareh}
\Theta_{\V}(\chi) := \{I_{\V}(\chi,\varphi) : \varphi\in
\Sc(\V(\A)) \} 
\end{equation}
is a subspace of ${\mathcal A}(\SLT(\A))$ where 
${\mathcal A}(\SLT(\A))$ is the space of genuine
automorphic forms on $\SLT(\A)$. Furthermore, the space $\Theta_{\V}(\chi)$ is
right invariant
under the action of $\SLT(\A)$ by right translation of the functions in
${\mathcal A}(\SLT(\A))$ \footnote{Technically, only the global Hecke
algebra is acting on function, so we are slightly abusing notation here}.
Our main result is a description of the space $\Theta_{\V}(\chi)$ for all
nontrivial automorphic characters $\chi$. This description will be given in terms
of another theta correspondence which we outline next.

Let $\U$ be the one dimensional space
with a quadratic form $Q(x)=mx^{2}$. The group $\GO(\U)$ is just $\{\pm 1\}$
and the adelic group $\GO(\U)(\A)=\prod_{v}\GO(\U)(\F_{v})$ is the
unrestricted direct product. The group $\GO(\U)(\A)$ is compact in the 
usual product topology  and $\GO(\U)(\F)$ is embedded diagonally in 
$\GO(\U)(\A),\;\; -1\mapsto (-1,-1,\ldots)$

A continuous character $\mu$ of $\GO(\U)(\A)$ determines a set of local
characters $\mu_{v}$ of $\GO(\U)(\F_{v})$
by restriction to the subgroup $\GO(\V)(\F_{v})\hookrightarrow \GO(\V)(\A)$
and continuity of $\mu$ implies that 
at almost all places $\mu_{v}={\bf 1}_{v}$, the trivial representation
of $\GO(\U)(\F_{v})$. At the places where
$\mu\neq {\bf 1}_{v}$ it
is the sign character $\textup{sgn}_{v}$. 
Let
\begin{equation} \label{S}
S=\{ v : \mu_{v}=\textup{sgn}_{v}\}.
\end{equation}
The automorphy condition on 
$\mu$ i.e. $\mu(-1)=1$ implies that $|S|$ is
even. Conversely, for any finite set of places $S$ where $|S|$ is even
there is a unique
automorphic character 
\[ 
\mu_{S}:=\left(\otimes_{v \in
S}\textup{sgn}_{v}\right)\otimes\left(\otimes_{v\notin S}{\bf 1}_{v}\right).
\]

The group $\SLT(\A)$ acts on the Schwartz 
space $\Sc(\U(\A))$ by the Weil
representation $\omega_{0}=\omega_{0,\psi}$ as in the case
above, the linear action of $\GO(\U)(\A)$ given by, 
$h_{0}\cdot\varphi_{0}(x)=\varphi(h_{0}^{-1}\cdot x)$ 
commutes with the action of
$\SLT(\A)$ on $\Sc(\U(\A))$. Using the commuting action, we define
for $h_{0}\in \GO(\U)(\A), \;\; g\in \SLT(\A),\;\; \varphi_{0}\in \Sc(\U(\A))$
the theta function
\[
\theta(g,h_{0},\varphi_{0}) =\sum_{x \in  \U(\F)}\omega_{0}(g)
(h_{0}\cdot \varphi_{0})(x)= 
\sum_{x \in \U(\F)}\omega_{0}(g)\varphi_{0}(h_{0}^{-1}\cdot x)
\]
As above, this function is left invariant by $\SL(\F)\times\GO(\U)(\F)$
and for an automorphic character $\mu$ of $\GO(\U)(\A)$,
we define the automorphic form on $\SLT(\A)$
\begin{equation}\label{simon}
 I_{\U}(\mu,\varphi_{0})(g)= 
\int\limits_{\GO(\U)(\F) \backslash \GO(\U)(\A)} 
\theta(g,h_{0},\varphi_{0}) \mu(h_{0})dh_{0}
\end{equation}
Note that in this case, $\GO(\U)(\A)$ is compact so the integral is
absolutely convergent.
The space
\begin{equation}\label{dan}
\Theta_{\U}(\mu)=\{ I_{\U}(\mu,\varphi_{0}) : \varphi_{0}\in\Sc(\U(\A))\}
\end{equation}
is a right $\SLT(\A)$ invariant subspace of ${\mathcal
A}(\SLT(\A))$. Such spaces for varying $\U$ and $\mu$, are 
colloquially referred to as ''coming from $\GO(1)$''.
We can now state our main result.
\begin{pf*}{Main Theorem}\label{main}
Let $\B$ be a quaternion algebra over a totally real field $\F$
Let $\chi$ be a quadratic character of $\A^{\times}/\F^{\times}$ given
by $\chi(x)=(x,-\kappa)_{\A}$, with $\kappa \in \F^{\times}$. 
Assume that $\chi \neq 1$ and that 
$S\neq \emptyset$ so that $\B$ is a division algebra.
Let $\Theta_{\V}(\chi)$ be the
corresponding automorphic representation of $\SLT(\A)$ constructed
above.
Then $\Theta_{\V}(\chi)$ is cuspidal irreducible and
\[
\Theta_{\V}(\chi)=\Theta_{\U}(\mu),
\]
where $\U$ is the one dimensional space with quadratic form
$Q(x)=\kappa x^{2}$ and $\mu=\mu_{S}$ is the automorphic
representation of $\GO(\U)(\A)$ determined by $S$.
\end{pf*}

\begin{rem}
The case of $\chi=1$ is covered by the classical Siegel-Weil formula. In that
case $\Theta_{\V}(1)$ is given by Eisenstein series whereas in our case 
$\Theta_{\V}(\chi)$ is a space of cusp forms. 
\end{rem}

One similarity between our
result and the Siegel Weil formula is that in the classical case, the
theta integral attached to $\varphi$ is identified with an explicit 
Eisenstein
series constructed from $\varphi$. Similarly, we are able to construct 
a function $\varphi_{0}\in \Sc(\U(\A))$ from $\varphi$ and identify
the theta lift attached to $\varphi$ with a theta lift attached to 
$\varphi_{0}$. This is made precise in the following
\begin{thm}\label{matching}
Suppose that $\B,\chi$ and corresponding $\U,\mu$ are as in the main theorem.
Then there is a unique map 
$\Sc(\V(\A))\rightarrow \Sc(\U(\A)),\; \varphi\mapsto \varphi_{0}$ 
such that
\[
I_{\V}(\chi,\varphi)(g)=I_{\U}(\mu,\varphi_{0})(g)
\]
for all $g \in \SLT(\A)$. More explicitly, for 
factorizable functions $\varphi=\otimes_{v} \varphi_{v}$ the map
$\varphi\mapsto \varphi_{0}$ is given by
$\varphi_{v}\mapsto\varphi_{0,v}$ where $\varphi_{0,v}$ is 
defined on $\F_{v}^{\times}$ by
\[
\varphi_{0,v}(r)=|r|_{v}\chi_{v}(r)
\int_{T_{x_{0}}(\F_{v})\backslash \HH(\F_{v})}
\varphi_{v}(r\, h^{-1}\cdot x_{0})\chi_{v}(\nu(h))dh.
\]
\end{thm}
We will say that the pair of functions $(\varphi,\varphi_{0})$ in
theorem \ref{matching} are matching. In fact
theorem \ref{main} follows from theorem \ref{matching} and the later
is proved via
a matching of Whittaker functions. 

It turns out that in our case, the automorphic form $I_{\V}(\chi,\varphi)$ is distinguished
\cite{GPS}, i.e., involves only one square class of Fourier
coefficients 
$\alpha$. That 
square class is the one where $\alpha=\kappa$, and
\[I_{\V}(\chi,\varphi)(n(b)g)
=\sum_{t\in \F^{\times}}
\psi^{\kappa t^{2}}(b)W_{\kappa}
(\begin{pmatrix}t & 0 \\ 0 & t^{-1} \end{pmatrix}g,\varphi).
\]
Note that although many authors (e.g. Gelbart and Piatetski Shapiro
\cite{GPS} for 
$\widetilde{{\rm GL}}_{2}(\A)$ and Vign\'eras, for classical
automorphic forms,) proved
that ``distinguished automorphic forms come from $\GO(1)$'',
their results do not apply in our case.

The automorphic forms in $\Theta_{\U}(\mu)$ are also distinguished.
For any $\varphi_{0}\in \Sc(\U(\A))$ 
the Fourier expansion of $I_{\U}(\mu,\varphi_{0})$ is
\[
I_{\U}(\mu,\varphi_{0})(n(b)g)
=
\sum_{t\in \F^{\times}}
\psi^{\kappa t^{2}}(b)W_{\kappa}(\begin{pmatrix}t & 0 \\ 0 & t^{-1} \end{pmatrix}g,\varphi_{0}),
\]
where, in fact, $W_{\kappa}(g,\varphi_{0})=\omega_{0}(g)\varphi_{0}(1)$.

To establish theorem \ref{matching} for a given $\varphi$ it suffices 
to construct $\varphi_{0}$ such that  $W_{\kappa}(g,\varphi_{0})= W_{\kappa}(g,\varphi)$
for all $g \in \SLT(\A)$. 
Thus the problem of matching reduces to a
problem of matching the Whittaker functions 
$
W_{\kappa}(g,\varphi)=W_{\kappa}(g,\varphi_{0}).
$
then, we are reduced to the problem of local matching:
given $\varphi_{v} \in \Sc(\V(\F_{v}))$, find $\varphi_{0,v}\in
\Sc(\U(\F_{v})$
such that 
\begin{equation} \label{lm}
W_{\kappa}(g_{v},\varphi_{v})=W_{\kappa}(g_{v},\varphi_{0,v})
\end{equation}
for all $g_{v}\in \SLT(\F_{v}))$.

The analysis of the local matching depends on the Bruhat decomposition
$\SLT(\F_{v})=\widetilde{B}\amalg \widetilde{B}w\widetilde{B}$ in the
following way. The
equality in 
equation \eqref{lm} for $g\in \widetilde{B}$ determines $\varphi_{0,v}$ uniquely for
$r \in \F_{v}^{\times}$ by the formula
\begin{equation}\label{vp}
\varphi_{0,v}(r)=|r|_{v}\chi_{v}(r)\int_{T_{x_{0}}(\F_{v})\backslash \HH(\F_{v})}
\varphi_{v}(r\, h^{-1}\cdot x_{0})\chi_{v}(\nu(h))dh.
\end{equation}
There are then two issues. First we must show that $\varphi_{0}$
defined on $\F_{v}^{\times}$ as in equation \eqref{vp}
can be extended to a Schwartz function on $\F_{v}$. Then, with this
Schwartz function 
$\varphi_{0,v}\in\Sc(\U(\F_{v}))=\Sc(\F_{v})$, we must show that
the equality in
\eqref{lm} holds on the big cell $\tilde{B}w\tilde{B}$ as well and hence on all of
$\SLT(\F_{v})$.

The calculation that establishes the local matching on the big cell amounts
to establishing
an identity involving orbital integrals, Gauss sums and Weil indexes.

For the global identity we must also check the fundamental lemma which
says that, in the unramified situation, the standard unramified vector
$\varphi_{v}^{0}\in \Sc(\V(\F_{v}))$ is matched by the standard unramified
vector  $\varphi_{0}^{0}\in \Sc(\U(\F_{v}))$.

One technical point which should be noted is that the usual theta
correspondence is between $\GO(\V)(\A)$ and $\SLT(\A)$ whereas we have
only discussed automorphic characters of ${\textup SO}(\V)(\A)$. However, in
our situation,
$\GO(\V)(\A)=\SO(\A)\times \mu\!\!\!\mu_{2}(\A)$ where
$\mu\!\!\!\mu_{2}=\{\pm 1\}$ and this product allows us to extend our
results
to $\GO(\V)(\A)$.
For an automorphic character $\eta$ of
$\mu\!\!\!\mu_{2}(\A)$, let
\[
\Theta_{\V}(\chi\otimes\eta)=\{ I_{\V}(\chi,\varphi):\varphi \in \Sc(\V(\A))^{\eta}\},
\]
where 
\[
\Sc(\V(\A))^{\eta}=\{\varphi \in \Sc(\V(\A)) :
\varphi(\delta x)=\eta(\delta)\varphi(x), \;\;\forall \;\;\delta \in \mu\!\!\!\mu_{2}(\A)\}.
\]
The space $\Theta_{\V}(\chi\otimes\eta)$ is the usual 
theta lift of the automorphic  character $\chi\otimes\eta$ of
$\GO(\V)(\A)$ and
as a consequence of our calculation we can determine it.
\begin{prop} \label{ue}
There is a unique character $\eta$ of $\mu\!\!\!\mu_{2}(\A)$ 
such that $\Theta_{\V}(\chi\otimes\eta)=\Theta_{\V}(\chi)$
and  $\Theta_{\V}(\chi\otimes\eta')=0$ for $\eta'\neq \eta$.
More precisely,
\[ 
\eta=\eta_{T}=\left(\otimes_{v \in T}\textup{sgn}_{v}\right)
\otimes\left(\otimes_{v\notin T}{\bf 1}_{v}\right),
\]
with
\[
T=\{ v : \chi_{v}(-1)\textup{Inv}(\B_{v})=-1 \}.
\]
\end{prop}

A final note about the type of representation we get is that
the local representations coming from $\GO(1)$ which we match are the
irreducible even or odd components of the Weil representation for a
one dimensional quadratic space. In fact, what we show is that, when
$\textup{Inv}(\B_{v})=1$, then in the local matching, $\varphi_{v}$ is
matched by an even $\varphi_{0}\in \Sc(\U(\F_{v}))$ and, when
$\textup{Inv}(\B_{v})=-1$, then $\varphi_{v}$ is
matched by an odd $\varphi_{0}\in \Sc(\U(\F_{v}))$. 
This means that the local component $\Theta_{\U}(\mu)_{v}$ of $\Theta_{\U}(\mu)$  
is the image of the irreducible representation of $\SLT(\A)$ on 
either the space of even Schwartz functions $\Sc(\U(\F_{v}))^{+}$
or the space of odd Schwartz functions $\Sc(\U(\F_{v}))^{-}$ according
to the invariant of $\B$. 
By theorem \ref{matching}
the local
component $\Theta_{\V}(\chi)_{v}$ of $\Theta_{\V}(\chi)$  coincides with
$\Theta_{\U}(\mu)_{v}$. Therefore the
local components $\Theta_{\V}(\chi)_{v}$ and the global representation
$\Theta_{\V}(\chi)$ are exactly the ones
coming from $\GO(1)$. 

In \cite{W2} Waldspurger considers the space 
${\mathcal A}_{00}(\SLT(\A))$ which
is the orthogonal complement in ${\mathcal A}_{0}(\SLT(\A))$
to the representation coming from $\GO(1)$ 
which he calls the
elementary Weil representations. These elementary representations are
exactly the ones that correspond to $\Theta_{\V}(\chi)$ in our case 
and so our result complements the classification carried out by Waldspurger.

Our result can be applied to the question of representations by spinor
generas to obtain results similar to those of R. Schulze-Pillot in
\cite{SP2}. Details will be included in a future article.

\subsection{Acknowledgments}
This paper is based on a thesis directed by Stephen Kudla who conceived
of the idea and directed me in carrying it out. I would like to take
this opportunity to thanks him for all his help.
Thanks are also due to
Jeff Adams for advice on the archemedian places and to Larry Washington for 
advice on Gauss and Kloosterman sums. All errors or misstatements that
remain are the responsibility of me alone.

\section{Notation} \label{dtl}
\subsection{Fields and characters}\label{ft}

When $\F_{v}$ is 
non archimedean we denote the ring of integers by 
$\OO_{v}$ and a generator of the maximal subring $\wp_{v}$ of $\OO_{v}$ by
$\varpi_{v}$. Let the modular degree of $\F_{v}$ over $\Q_{p}$ be $r_{v}$ so
that $|\varpi_{v}|_{v}=q=p^{r_{v}}$

Fix $\psi$-a non trivial additive character of $\A=\A_{\F}$ which is
trivial on $\F$. The local
components of $\psi$, will be denoted by $\psi_{v}$ except when it is 
clear from context that $\psi$ means its local component.

The character $\psi$ is defined in terms of the following 
parameters. We can define a map $\lambda$ as the 
composition of the maps
$\F_{v} \rightarrow \Q_{p}\rightarrow \Q_{p}/\Z_{p} \hookrightarrow \Q/\Z$.
where the map from $\F_{v}$ to $\Q_{p}$ is the trace map.
We then define characters of $\F_{v}$ 
\[
\psi_{0,v}(x)=\exp(-2\pi i \lambda(x))
\]
The product of the local characters $\psi_{0}=\prod_{v}\psi_{0,v}$ 
is a character of $\A/F$ and
our fixed character $\psi$ and its local factors are defined as 
\begin{equation}\label{defpsi}
\psi(x)=\psi_{0}(\epsilon x),\;\; \;\;\;
\psi_{v}(x)=\psi_{0,v}(\epsilon x) \;\;\;
\end{equation}
for some $\epsilon\in \F^{\times}$. For $\alpha \in \F^{\times}$, let
$\psi^{\alpha}(x)=\psi(\alpha x )$. Then $\{ \psi^{\alpha}\; :\;
\alpha \in \F^{\times}\}$ is the space of all non trivial characters
on $\A/\F$.

The conductor of $\psi_{0,v}$ is the smallest power of $\wp_{v}$ on
which $\psi_{0,v}$ is trivial. It is equal to
the inverse
different ${\mathcal D}^{-1} = \{ x \in \F_{v}\; :\; \textup{Tr}(x) \in
\Z_{p}\}$. The differental exponent is $d$, it is defined such that
${\mathcal D}=\wp^{-d}$. Note that the conductor of $\psi$ is 
$\epsilon^{-1}{\mathcal D}^{-1}=\epsilon^{-1}\wp^{-d}$.

\subsection{The quaternions} \label{quat}
$\B$ is a quaternion algebra over $\F$, and 
$\B_{v}=\B\otimes_{\F}\F_{v}$. The Invariant $\textup{Inv}(\B_{v})=+1$ if 
$\B_{v}\simeq \textup{M}_{2}(\F_{v})$ and $\textup{Inv}(\B_{v})=-1$ if
$\B_{v}$ is the (unique up to isomorphism) division quaternion algebra 
over $\F_{v}$.

The isomorphism class of $\B$ is uniquely determined by the set 
of places $S=S_{\B}:=\{v: \textup{Inv}(\B_{v})=-1\}$.

\begin{defn}
$\V$ is the 3-dimensional $\F$-vector space
\[ 
\V=\{x \in \B : \textup{Tr}(x)=0 \}
\]

We also define $\V_{v}=\V\otimes_{\F}\F_{v}$ and 
$\V(\A)=\V\otimes_{\F}\A$. We will often abuse notation and
write $\V$ instead of $\V_{v}$ when the meaning of $\V$ is clear from
context.
\end{defn}

\begin{defn}
The group $\HH$ is the group of invertible quaternions
\[ 
\HH=\{x \in \B : \nu(x)\neq 0 \}
\]
\end{defn}
\begin{defn}
The group $\T_{x}\subset\HH$ is the centralizer of $x$
\[ 
\T_{x}=\{y \in \HH : xy=yx \}
\]
\end{defn}

The standard definitions of Schwartz spaces are
\begin{defn} \label{scdef}
When $\V$ is defined over a non archimedean field, 
the Schwartz space $\Sc(\V)$ is the set of locally constant compactly
supported functions on $\V$.

 When $\V$ is defined over an archimedean field, 
the Schwartz space $\Sc(\V)$ is the set of infinitely differentiable
rapidly decreasing functions. When $\V$ is split 
and archimedean, the ``restricted Schwartz space'' is the 
subspace of functions of the form 
\[
f\begin{pmatrix} u & v \\ w & -u \end{pmatrix}=
P(u,v,w)\exp(-2\pi |\epsilon|(u^{2}+\frac{v^{2}}{2}+\frac{w^{2}}{2})
\]
where $P$ is a polynomial.
 
When $\V$ is non split and archimedean, 
the ``restricted Schwartz space'' is the space of functions
\[
f(x)=
P(x)\exp(-2\pi |\epsilon|\nu(x))
\]
Where $P$ is a polynomial and $x=(u,v,w)\in \V$.
\end{defn}
This section concludes with
some basic geometric properties of $\B$.

\begin{lem}\label{comp}
If $\B_{v} =\M_{2}(\F_{v})$ then 
$\ZZ(\F_{v}) \backslash \HH(\F_{v})$ is not compact
and if 
$\B_{v}={\mathbb H}(\F_{v})$ then 
$\ZZ(\F_{v}) \backslash \HH(\F_{v})$ is compact.
\end{lem}

\begin{lem} The quotient
$\ZZ(\A) \backslash \HH(\A)$ 
is compact unless $S=\emptyset$.
\end{lem}
\begin{pf}
This follows from corollary 3.2 in \cite{VI}. The product of the reduced norms
is isotropic if and only if the reduced norm is isotropic at every $v$.
\end{pf}

\begin{lem} \label{ext}
For $x \in \V $ we have
$\T_{x}=(\F(x))^{\times}$
\end{lem}

\begin{rem}
Lemma \ref{ext} also holds over all local fields $\F_{v}$ and over the
ring of Adeles, $\A$.
\end{rem}

\begin{lem} \label{compactness}
 If $\B=B_{v}$ then  $\nu(x) \notin -(\F^{\times})^{2}$
if and only if  $\ZZ \backslash \T_{x}\simeq \ZZ \backslash \F(x)$  is compact.
\end{lem}

\subsection{$\SLT(\A)$}
The group $\SLT(\A)$ is the unique non trivial 2-fold topological
central extension of $\SL_{2}(\A)$
$$ 1 \rightarrow \{\pm\} \rightarrow \SLT(\A)
\rightarrow \SL_{2}(\A) \rightarrow 1
$$
Group elements are realized as pairs 
$(g,\epsilon)\;\; g\in \SL_{2}(\A) \;\; \epsilon \in  \{\pm 1\}
$. Multiplication is defined by:
\[
(g_{1},\epsilon_{1})(g_{2},\epsilon_{2})=
(g_{1}g_{2},\epsilon_{1}\epsilon_{2}\beta(g_{1},g_{2}))
\] 
Where $\beta$ is a 2-cocycle on $\SL_{2}(\A)$.

The generators of  $\SLT(\A)$ are:
\begin{equation}\label{gendef}
m(a)=\left(\begin{array}{cc}
        a & 0 \\
        0 & a^{-1}
\end{array}\right)\;\; 
n(b)= \left(\begin{array}{cc}
        1 & b \\
        0 & 1
\end{array}\right) \; \;
 w = \left(\begin{array}{cc}
        0 & 1 \\
        -1 & 0
\end{array}\right) \;\; (1,-1)
\end{equation}
When calculating the action of group elements by the Weil
representation (see below) we will drop the $\epsilon$ from the notation.

\subsection{$\GO(V)(\A)$ and $\SLT(\A)$} 
Recall that we have fixed a division quaternion algebra $\B$ with
associated quadratic space $\V$ and hence indentified 
$\HH=\B^{\times}=\textup{GSpin}(\V)$.

The actions of $\SLT(\A)$ and $\HH(\A)$ on $\Sc(V(\A))$ are the following: 
$\HH(\A)$ acts on
elements of $\V(\A)$ by 
conjugation and on elements of $\Sc(V(\A))$ by 
acting on their arguments i.e. for
$\varphi\in \Sc(\V(\A)), \;\;\;\; x\in \V(\A)$, and $ h \in \HH(\A)$,
the action is
\begin{gather*}
h\cdot \varphi(x) = \varphi(h^{-1}\cdot x)= \varphi(h^{-1} x h)
\end{gather*}
 
The Weil representation of $\SLT(\A)$ as described in 
Kudla's notes (\cite{KU} p. 32-34)
is determined by an additive
character $\psi$ of $\A$ which we fixed in chapter \ref{dtl}
and by the quadratic space $(V(\A),\nu)$, determined by the algebra 
$\B(\A)$.

The formulas for the actions of the
generators of $\SLT(\A)$ on Schwartz functions are:
\begin{eqnarray}
\omega(n(b))\varphi(x)&=&\psi(b\nu(x))\varphi(x) \label{wnb}\\
\omega(m(a))\varphi(x)&=& |a|^{\frac{3}{2}}
(a,-\det(\V))\gamma(a,\psi^{\frac{1}{2}})^{-1}\varphi(ax)\label{wma}\\
\omega((w,\epsilon))\varphi(x)&=& \epsilon
\gamma(\psi\circ \V)^{-1}\int_{\V(\A)} \varphi(y)\psi(-(x,y))dy. \label{ww}
\end{eqnarray}

Here $\gamma(a,\psi^{\frac{1}{2}})$ is the quotient of the Weil
indices
$\frac{\gamma(\psi^{\frac{a}{2}})}{\gamma(\psi^{\frac{1}{2}})}$. See 
the appendix in 
\cite{RA} for the definitions and properties of 
the Weil index. Also, $\det(\V)$ is the determinant of the matrix of
the quadratic form on $\V$ and 
\[
\gamma(\psi\circ \V)=
\gamma(\det(\V),\psi^{\frac{1}{2}})\gamma(\psi^{\frac{1}{2}})^{3}\epsilon(\V)
\]
where $\epsilon(\V)=\pm 1$ is the Hasse invariant of $\V$, see \cite{SE}. 
The measure $dy$ in equation \eqref{ww} is the unique Haar measure
which makes the transform $\hat{f}(x)=\int_{\V}f(y)\psi(-(x,y))dy$
self dual. For more details see \cite{KU}.

\subsection{$\GO(\U)(\A)$ and $\SLT(\A)$}
For the 1-dimensional quadratic space $\U=\F$ with $Q(x)=mx^{2},\;
m\in \F^{\times}$
the Weil representation
$\omega_{0}$ of $\SLT(\A)$, commutes with the action of
$\GO(\U)(\A)$.This representation is given by the formulas
\begin{eqnarray}
\omega_{0}(n(b))\varphi_{0}(x)&=&\psi(bmx^{2})\varphi_{0}(x) \label{w0nb}\\
\omega_{0}(m(a))\varphi_{0}(x)&=&
|a|^{\frac{1}{2}}(a,2m)\gamma(a,\psi^{\frac{1}{2}})^{-1}\varphi_{0}(ax)\label{w0ma}\\
\omega_{0}((w,\epsilon))\varphi(x)&=& \epsilon
\gamma(\psi\circ \U)^{-1}\int_{\U(\A)} \varphi_{0}(y)\psi(-2mxy)dy \label{w0w}
\end{eqnarray}
As above  $\gamma(a,\psi^{\frac{1}{2}})$ is the quotient of the Weil
indices
$\frac{\gamma(\psi^{\frac{a}{2}})}{\gamma(\psi^{\frac{1}{2}})}$
and 
\[
\gamma(\psi\circ \U)
=\gamma(\det(\U),\psi^{\frac{1}{2}})\gamma(\psi^{\frac{1}{2}})\epsilon(\U)
\]

\section{Fourier expansions and Whittaker models}

Our matching of the spaces $\Theta_{\V}(\chi)$ and $\Theta_{\U}(\mu)$
will be reduced to matching their 
 Whittaker models so we first realize the two spaces as such.

\subsection{Whittaker model of $\Theta_{\V}(\chi)$}
We will now calculate the Fourier expansion of functions 
$I_{\V}(\chi,\varphi)\in \Theta_{\V}(\chi)$. This calculation will show that,
like $\Theta_{\U}(\mu)$, the space $\Theta_{\V}(\chi)$ is distinguished.
We start by writing the automorphic forms as sums.
\begin{align}
I_{\V}(\chi,\varphi)(g)&=\int\limits_{\HH(\F)\ZZ(\A)\backslash \HH(\A)}
\theta(g,h,\varphi)\chi(\nu(h))dh \nonumber \\
&=\int\limits_{\HH(\F)\ZZ(\A)\backslash \HH(\A)}
\sum_{x \in \V(\F)}\omega(g)\varphi(h^{-1}\cdot x) \chi(\nu(h))dh \nonumber \\
&= \int\limits_{\HH(\F)\ZZ(\A)\backslash \HH(\A)}
\sum\limits_{\substack{\delta \in \F\\ \delta \neq -\Box^{\times}} }
\sum\limits_{\substack{x\in \V(\F)\\ \nu(x)=\delta}}
\omega(g)\varphi(h^{-1}\cdot x) \chi(\nu(h))dh \label{dl}
\end{align}
We make the convention that $\Box$ should be read as 
``a square'' and $\Box^{\times}$ be read as ``a non zero square``. 
As mentioned in sections \ref{quat}, we exclude the case 
of the unramified global algebra $\B(\A)=\M(2,\A)$ and therefore
by Lemmas \ref{comp} and \ref{compactness} the quotient
$\HH(\F)\ZZ(\A)\backslash \HH(\A)$ is compact and the integral in
equation \eqref{dl} is absolutely convergent. Also, $\V$ is anisotropic
so that
the norms of $x \in \V$ are not equal to a
negative of a non zero square.
Therefore $-\delta$ is not equal to a non zero square.
\begin{lem}
The term $\delta=0$ in the sum in equation \eqref{dl} is zero.
\end{lem}
\begin{pf}
The space $\V$ is
anisotropic so the only point of norm 0 is 0 and the the $\delta=0$-term
in the sum is the integral of the character $\chi\circ\nu$ on the 
compact space $\HH(\F)\ZZ(\A)\backslash \HH(\A)$ so the term is zero. 
\end{pf}
The integral in equation \eqref{dl} 
can be switched with the summation since the domain of 
integration is compact and the function has compact support.
Thus we have
\[
I_{\V}(\chi,\varphi)(g)= 
\sum\limits_{\substack{\delta \in \F^{\times}\\\delta \neq - \Box}} 
\; \int\limits_{\HH(\F)\ZZ(\A)\backslash \HH(\A)}
 \sum\limits_{\gamma \in \T_{x_{\delta}}(\F)\backslash \HH(\F)}
\!\!\!\!\!\!\omega(g)\varphi(h^{-1}\gamma^{-1}\cdot
x_{\delta}) \chi(\nu(h))dh,
\]
where $x_{\delta}\in \V(\F)$ is some point such that $\nu(x_{\delta})=\delta$
Notice that the summation over the set $\{ x\in \V(\F)\;:\; \nu(x)=\delta\}$ 
is replaced by the sum over $\T_{x_{\delta}}(\F)\backslash \HH(\F)$.
 Since $\B$ is a division algebra
and $\V$ is anisotropic $\HH(\F)$ acts transitively on
$\{ x \in \F \,|\, \nu(x)=\delta\}$ and so the choice of 
$x_{\delta}$ does not matter. 
Also, since $\chi$ is
trivial on the rational points $\chi(\nu(\gamma h)) =\chi(\nu(h))$ for
$\gamma \in \HH(\F)$ so
we can include $\gamma$ in the argument of $\chi\circ \nu$.
\[
I_{\V}(\chi,\varphi)(g)=
\!\sum\limits_{\substack{\delta \in \F^{\times} \\ \delta \neq - \Box}} 
\; \int\limits_{\HH(\F)\ZZ(\A)\backslash \HH(\A)}
 \sum\limits_{\gamma \in \T_{x_{\delta}}(\F)\backslash \HH(\F)}
\!\!\!\!\!\!\!\omega(g)\varphi(h^{-1}\gamma^{-1}\cdot x_{\delta}) \chi(\nu(\gamma h))dh 
\]
We now combine the inner sum and integral by a change of variables
$h_{0}=\gamma h$.
\[
I_{\V}(\chi,\varphi)(g)=
\sum_{\substack{\delta \in \F^{\times}\\ \delta \neq - \Box}}
\int_{\T_{x_{\delta}}(\F)\ZZ(\A)\backslash \HH(\A)} 
\omega(g)\varphi(h_{0}^{-1}\cdot
x_{\delta}) \chi(\nu(h_{0}))dh_{0} 
\]
The next step is to factor $h_{0}=h_{1}q_{1}$ with
$h_{1}\in A_{\delta}:=\T_{x_{\delta}}(\F)\ZZ(\A)\backslash \T_{x_{\delta}}(\A)$ 
and, 
$q_{1} \in B_{\delta}:=\T_{x_{\delta}}(\A)\backslash \HH(\A) $.
This gives, 
\begin{align}
I_{\V}(\chi,\varphi)(g)
&=\sum_{\substack{\delta \in \F^{\times}\\ \delta \neq -\Box }}
\int_{B_{\delta}}
\int_{A_{\delta}}
\omega(g)\varphi(q_{1}^{-1}h_{1}^{-1}\cdot
x_{\delta}) \chi(\nu(h_{1}q_{1}))dh_{1}\,dq_{1} \\
&=\sum_{\substack{\delta \in \F^{\times}\\ \delta \neq -\Box }}
\int_{B_{\delta}}
\omega(g)\varphi(q_{1}\cdot x_{\delta}) \chi(\nu(q_{1}^{-1}))dq_{1}
\int_{A_{\delta}}\chi(\nu(h_{1}))dh_{1} \label{twoint}
\end{align}
To proceed, we need the
following description of
$\T_{x_{\delta}}(\A)$.
\begin{lem} \label{ker}
Let $\chi$ is a nontrivial quadratic character of $\F\backslash \A$ 
given by the Hilbert symbol $\chi(r)=(r,-\kappa)_{\A}$ and
$\chi_{v}(r)=(r,-\kappa)_{\F_{v}}$. Let $\delta =\nu(x_{\delta})\neq 0$ with 
$x_{\delta}\in \V(\F)$.
\begin{enumerate}
\item 
$\T_{x_{\delta}}(\F_{v})\subset
\textup{ker}(\chi_{v}\circ\nu_{v})$
if and only if $\delta \in \kappa(\F_{v}^{\times})^{2}$
\item
$\T_{x_{\delta}}(\A)\subset \textup{ker}(\chi\circ\nu)$ 
if and only if $\delta \in \kappa (\F^{\times})^{2}$
\end{enumerate}
\end{lem}
\begin{pf}
The first result is contained in Proposition 4.8.2 in \cite{B} when we
notice that
Since $\chi_{v}(a)=(a,-\kappa)_{\F_{v}}$ 
the kernel of $\chi_{v}$ is the set 
\[ \textup{ker}(\chi_{v})= \nu(\F_{v}(\sqrt{-\kappa})^{\times})
\]
and since 
$\T_{x_{\delta}}(\F_{v})=(\F_{v}(x_{\delta}))^{\times}=
(\F_{v}(\sqrt{\delta})^{\times}$,
\[
\nu(\T_{x_{\delta}}(\F_{v}))=\textup{ker}(\chi_{\delta})
\]
where $\chi_{\delta}(a)=(a,-\delta)_{\F_{v}}$. 

The second statement follows from the fact that if ${\mathcal K}_{1}$
and ${\mathcal K}_{2}$ are quadratic extensions of $\F$ then they
determine characters $\chi_{1}$ and $\chi_{2}$, where 
$\ker(\chi_{i})=\nu({\mathcal K}_{i}(\A)^{\times})$, and 
$\ker(\chi_{1})=\ker(\chi_{2})$ if and only if 
${\mathcal K}_{1}\simeq {\mathcal K}_{2}$.
\end{pf}
\begin{cor}\label{ramtriv}
If $\chi_{v}=1$ for some place $v\in S$, the set of ramified places for $\B$
then $\Theta_{\V}(\chi)=0$.
\end{cor}
\begin{pf}
If $\chi_{v}$ is trivial at a ramified place $v$, the square class of 
$\kappa$ is $-(\F_{v}^{\times})^{2}$. But since $\delta \neq -\Box $ 
and $\delta$ is not a norm of $\V(\F_{v})$,
$\delta$  is not in the
square class of $\kappa$.
Therefore, by Lemma \ref{ker} 
$\T_{x_{\delta}}(\A)\not\subset \textup{ker}(\chi\circ\nu)$
and the integral over $A_{\delta}$ in equation \eqref{twoint} will 
be zero for any $g$ and $\varphi$.
\end{pf}

Using Lemma \ref{ker} we can see that the inner integral in 
\eqref{twoint} is nonzero exactly when $\delta \in \kappa (\F^{\times})^{2}$.
If we pick $x_{0}$ such that $\nu(x_{0})=\kappa$ then
elements with norm $\alpha^{2}\kappa$ are represented by $\alpha x_{0}$, therefore
\begin{equation}  \label{I_chi}
I_{\V}(\chi,\varphi)(g)=
\sum\limits_{\alpha \in \F^{\times} }
\int\limits_{\T_{x_{0}}(\A)\backslash \HH(\A)}\omega(g)
\varphi(q^{-1}\cdot \alpha x_{0})\chi(\nu(q))dq 
\end{equation}
where we normalized the volume such that 
$$\textup{Vol}(\T_{x_{0}}(\F)\ZZ(\A)\backslash
\T_{x_{0}}(\A))=\textup{Vol}(\T_{\delta x_{0}}(\F)\ZZ(\A)\backslash
\T_{\delta x_{0}}(\A))=1$$
When we replace $g$ in equation \eqref{I_chi} by $n(b)g$ then the 
$\alpha$ term in the sum is multiplied by $\psi(\alpha^{2}\kappa b)$. 
Therefore
the sum in \eqref{I_chi} is actually the Fourier expansion of 
$I_{\V}(\chi,\varphi)$. 
We make this explicit in the following
definitions:
\begin{defn}\label{wvarphi}
The Whittaker function corresponding to
the form  $I_{\V}(\chi,\varphi)$ is the function
\[
W_{\varphi}(g)=\int\limits_{\T_{x_{0}}(\A)\backslash \HH(\A)}\omega(g)
\varphi(q^{-1}\cdot x_{0})\chi(\nu(q))dq.
\]
where $x_{0}\in \V_{\F}$ such that $\nu(x_{0})=\kappa$
\end{defn}

\begin{defn} \label{ed}
The Whittaker model of the automorphic representation $\Theta_{\V}(\chi)$
is the space 
\[
W_{\V}(\chi):=\{W_{\varphi}(g) : \varphi \in \Sc(\V(\A)) \} 
\]
\end{defn}
We collect the result of the above discussion in the proposition
\begin{prop}\label{dist}
Let $\chi\neq 1$ then
\begin{enumerate}
\item The function $I_{\V}(\chi,\varphi)$ is a
cuspidal automorphic form and therefore 
\[
\Theta_{\V}(\chi)\subset {\mathcal A}_{0}(\SLT(\A)).
\]

\item  The function $I_{\V}(\chi,\varphi)$ is distinguished in the
  sense of Gelbart and Piatetski-Shapiro \cite{GPS} with
Fourier coefficients in the square class of $\kappa(\F^{\times})^{2}$.
\end{enumerate}
\end{prop}
The information we have so far, allows us to determine some
information about the space $\U$.
\begin{rem}\label{setm}
In order that $\Theta_{\V}(\chi)\subset \Theta_{\U}(\mu)$ it is necessary
that the quadratic form on $\U$ is defined by $Q(x)=\kappa x^{2}$.
\end{rem}

\subsection{The Whittaker model of $\Theta_{\U}(\mu)$}
In this section we compute the Fourier expansion of the function
$I_{\U}(\chi,\varphi_{0})$.
First, we write 
\begin{align*}
I_{\U}(\chi,\varphi_{0}) &= \int_{\GO(\U)(\F) \backslash \GO(\U)(\A)} \sum_{x \in \F}
\omega_{0}(g)\varphi_{0}(h_{0}^{-1}\cdot x)\mu(h_{0})dh_{0}\\
&= I_{\U}(\mu,\varphi_{0})(g)=\frac{1}{2}\sum_{x \in \F}
\int_{\GO(\U)(\A)} 
\omega_{0}(g)\varphi_{0}(u^{-1}\cdot x)\mu(u)du 
\end{align*}

The group $\GO(\U)(\A)$ is compact, so the representation on 
$\Sc(\V(\A))$ is a direct sum of isotypic components and it is not hard
to see that the inner integral is the volume of $\GO(\U)(\A)$ times the projection 
from $\Sc(\U(\A))$ onto the $\mu=\mu_{S}$ isotypic subspace. Recall that
at places $v\in S$ the
factor $\mu_{v}$ is the sign representation and at places $v \notin S$
the factor $\mu_{v}$ is the trivial character. Therefore the 
$\mu=\mu_{S}$ isotypic subspace is the subspace
\begin{equation} \label{parity}
\Sc(\U(\A))_{\mu}=(\otimes_{v\in S}\Sc(\U(\F_{v}))^{-})\otimes
(\otimes_{v\notin S}\Sc(\U(\F_{v}))^{+}),
\end{equation}
where
$\Sc(\U(\F_{v}))^{-},\;\Sc(\U(\F_{v}))^{+}$ are the subspaces of odd and even
functions respectively.

This means that we can replace $\varphi_{0}$ by its projection to
$\Sc(\U(\A))_{\mu}$ which we also denote by $\varphi_{0}$ and then
\begin{align*}
I_{\U}(\mu,\varphi_{0})(g)
&=\frac{1}{2}\sum_{x \in \F} \omega_{0}(g)\varphi_{0}(x).
\end{align*}
The space $\Theta_{\U}(\mu)$ is isomorphic to 
$\Sc(\U(\A))_{\mu}$.

 To construct the Whittaker model for $\Theta_{\U}(\mu)$ we
consider the Fourier expansion of the forms $I_{\U}(\mu,\varphi_{0})$.
Let $g\in \SLT(\A)$ be any element then for any $b \in \A$
\begin{align*}
I_{\U}(\mu,\varphi_{0})(n(b)g)
&= \frac{1}{2}\sum_{x \in \F}\omega_{0}(n(b)g)\varphi_{0}(x)\\
&= \frac{1}{2}\sum_{x \in \F}\psi(bmx^{2})\omega_{0}(g)\varphi_{0}(x)
\end{align*}
The term $x=0$ can be removed from the sum since
$\varphi_{0}\in\Sc(\U(\A))_{\mu}$ and $\mu \neq 1 $, so $\varphi_{0}(0)=0$. 

By equation \eqref{w0ma} $\omega_{0}(m(t)g)\varphi_{0}(1)=\omega_{0}(g)(t)$
so we can rewrite the
Fourier expansion as
\begin{align*}
I_{\U}(\mu,\varphi_{0})(n(b)g)=\frac{1}{2}
\sum_{x \in \F^{\times}}\psi^{mx^{2}}(b)\omega_{0}(m(x)g)\varphi_{0}(1)
\end{align*}
Thus $I_{\U}(\mu,\varphi_{0})$ is cuspidal and the
Fourier expansion is completely determined by the function
$\omega_{0}(g)\varphi_{0}(1)$.
Notice that the character $\psi^{mx^{2}}$ appears twice in the sum, at $\pm
x$ and that $\varphi_{0}\in \Sc(\V(\A))_{\mu}$ is even so the coefficient of the character
is actually $\omega_{0}(m(x)g)\varphi_{0}(1)$.
The above calculation justifies the definitions:
\begin{defn} \label{wvarphi0}
The Whittaker function corresponding to the automorphic form $I_{\U}(\mu,\varphi_{0})$
is
\[ 
W_{\varphi_{0}}(g)=\omega(g)\varphi_{0}(1)
\]
\end{defn}

\begin{defn} \label{whitt}
The global Whittaker model of  $\Theta_{\U}(\mu)$
is the space 
\[
W_{\U}(\mu):= \{W_{\varphi_{0}} : \varphi_{0} \in \Sc(\U(\A))_{\mu} \} 
\]
\end{defn}
Recall that $m$ defined the
quadratic form on $\U,\; Q(x)=mx^{2}$.
We note a result which was evident in the Fourier expansion of $I_{\U}(\mu,\varphi_{0})$.
\begin{lem} \label{dist1}
The automorphic forms $I_{\U}(\mu,\varphi_{0})$ are cuspidal and 
distinguished in the sense of  Gelbart and Piatetski-Shapiro
\cite{GPS}. The characters in the Fourier expansion of 
$I_{\U}(\mu,\varphi_{0})$ are all in the square class of $m$.
\end{lem}

\subsection{Factorizability}\label{factorizability}

If $\varphi\in \Sc(\V(\A))$ is factorizable, 
we write $\varphi = \otimes_{v}\varphi_{v}$ and then
\begin{equation*}
W_{\varphi}(g)
=\int\limits_{\T_{x_{0}}(\A)\backslash \HH(\A)}\omega(g)
\varphi(q^{-1}\cdot x_{0})\chi(\nu(q))dq
=\prod_{v}W_{\varphi_{v}}(g_{v})
\end{equation*}
where 
$
W_{\varphi_{v}}(g_{v})=\int\limits_{\T_{x_{0}}(\F_{v})\backslash \HH(\F_{v})}\omega_{v}(g_{v})
\varphi(q_{v}^{-1}\cdot x_{0})\chi_{v}(\nu_{v}(q_{v}))dq_{v}
$
and the $\psi^{\kappa}$-Whittaker space of all such functions is
\[
W_{\psi^{\kappa}}(\V,\chi) := \{W_{\varphi_{v}} \;:\; \varphi \in \Sc(\V(\F_{v}))\}. 
\]
Returning to the case of one dimensional
Schwartz functions $\varphi_{0}$ we have a similar factorization where
we write an arbitarary function as a sum of factorizable functions.
 if $\varphi_{0}\in \Sc(\U(\A))$ is factorizable we write
$\varphi_{0} = \prod_{v}\varphi_{0,v}$ and then
\[
W_{\varphi_{0}}(g)=\omega_{0}(g)\varphi_{0}(1)
=\prod_{v}W_{\varphi_{0,v}}(g_{v})
\]
where
\[
W_{\varphi_{0}}(g_{v})=\omega_{0,v}(g_{v})\varphi_{0,v}(1)
\]      
and the $\psi^{\kappa}$-Whittaker space (by Lemma \ref{setm}) of all such functions is
\[
W_{\psi^{\kappa}}(\U,\mu) := \{W_{\varphi_{0,v}} \;:\; \varphi \in \Sc(\U(\F_{v}))\}. 
\]

From here on, our calculation will consist of local
matching. Therefore we will drop the $v$ subscript from most
notation. Unless stated otherwise, all formulas should be
understood to be over a completion $\F_{v}$. We will also write 
$T_{0}$ instead of $\T_{x_{0}}$.

\section{Local Matching}

\subsection{Structural preliminaries}

The various lemmas in this section will be used in the later sections.

\begin{rem}
By Lemma \ref{ramtriv} we can assume that, whenever $\chi=1$, then 
$\B_{v}=\M_{2}(\F_{v})$.
\end{rem}

\begin{lem}\label{conji}
With notation as above, for every place $v$ there is an element $\xi \in \HH$ 
such that $\xi^{-1}\cdot x_{0}=-x_{0}$ and $\chi(-\nu(\xi))=\textup{Inv}(\B_{v})$.
\end{lem}
\begin{pf}
By equation (1) on page 1 and corollary 2.2 on page 6. in \cite{VI} we
can realize the quaternion
algebra as $\F_{v}(x_{0}) + \F_{v}(x_{0})\xi$ with $\xi$ such that if 
$m \in \F_{v}(x_{0})$ then $\xi m=\bar{m}\xi$. Where $\bar{m}$ is the
conjugate of $m$ in $\F_{v}(x_{0})$. It follows that
$\xi x_{0}\xi^{-1}=\bar{x}_{0}\xi \xi^{-1}=\bar{x}_{0}=-x_{0}$.
Also, by corollary 2.4 on page 6 in \cite{VI} we have
$\textup{Inv}(\B)=1$ if and only if 
$-\nu(\xi) \in \nu(\F_{v}(x_{0}))=\textup{ker}(\chi)$. This means that 
$\textup{Inv}(\B)=\chi(-\nu(\xi))$.
\end{pf}

\begin{defn} \label{halpha}
For $\alpha \in \F$, let $h_{\alpha}\in \HH$ be any element such that 
$h_{\alpha}^{-1}x_{0}=\alpha x_{0} + x$ with $x_{0}\perp x$.
\end{defn}
\begin{lem} \label{cosets}
Let $\chi\neq 1$, then the coset space 
$\T_{0}\backslash \HH$ is equal to
the disjoint union of right cosets $h_{\alpha}\T_{0}$ over the set
of $\alpha\in\F^{\times}$ such that 
$$
\chi(\alpha^{2}-1)=\textup{Inv}(\B) \;\;\textup{or} \;\;
\alpha=\pm 1.
$$
\end{lem}
\begin{pf}
If $h$ and $h'$ are such that for some $\alpha$, 
$h^{-1}x_{0}=\alpha x_{0} + x$ and
$h'^{-1}x_{0}=\alpha x_{0} + x'$ with $x$ and $x'$ perpendicular to
$x_{0}$ then $\nu(x)=\nu(x')$ and there exists an element of
the orthogonal group of the orthogonal complement to the span of $x_{0}$ 
taking $x$ to $x'$. This means that
$h\T_{x_{o}}=h'\T_{x_{0}}$ so the decomposition into   
right cosets indexed by $\alpha$ is proved.

To verify the second part of the claim realize $\B_{v}$ as 
$\F_{v}(x_{0}) + \F_{v}(x_{0})\xi$ as in the proof of \ref{conji}

To carry out the
calculation we use the orthogonal basis
$\{x_{0},\xi,x_{0}\xi\}$ for $\V$ and
write 
an arbitrary element $x$ in terms of this basis:
$x=\alpha x_{0}+\beta \xi +\gamma x_{0}\xi $. In terms of those
coordinates,
$\nu(x)=\alpha^{2}\kappa +\beta^{2}\nu(\xi) +\gamma^{2}\kappa\nu(\xi)$.
If $x$ is in the orbit of $x_{0}$ then
\[
\nu(x)=\alpha^{2}\kappa +\beta^{2}\nu(\xi) +\gamma^{2}\kappa\nu(\xi)=\kappa
\]
This means that $(\alpha^{2}-1)\kappa=-\nu(\xi)(\beta^{2}+\gamma^{2}\kappa)$
and then $\chi(\alpha^{2}-1)\chi(\kappa)=\chi(-\nu(\xi))\chi(\beta^{2} +\gamma^{2}\kappa)$.
Now, since $\chi(\beta^{2}+\gamma^{2}\kappa)=\chi(\kappa)=1$ we are done.
\end{pf}

The second lemma describes the values of $\chi(h_{\alpha})$:
\begin{lem} \label{chihalpha}
If $h_{\alpha}$ is defined as in \ref{halpha} then 
\begin{equation}
\chi(\nu(h_{\alpha}))=\chi(2(\alpha +1))
\end{equation}
\end{lem}
\begin{pf}
we extend the basis defined in the proof of Lemma \ref{cosets} for
$\V$ to all of $\B_{v}$ and write 
$h_{\alpha}=a +bx_{0}+c\xi +dx_{0}\xi$ then by multiplying out the 
conjugation action of $h_{\alpha}$:
and collecting the coefficients of 
$x_{0}$ we get
\[
\alpha = \frac{a^{2}+b^{2}\kappa - c^{2}\nu(\xi)-d^{2}\kappa \nu(\xi)}
        {a^{2}+b^{2}\kappa +c^{2}\nu(\xi) +d^{2}\kappa \nu(\xi)}
\]
From that we get that $\alpha + 1 = 
\frac{2(a^{2} +b^{2}\kappa)}{\nu(h_{\alpha})}$
and since $\chi(a^{2} +b^{2}\kappa)=1$ we get
\[\chi(2(\alpha +1))=\chi(\nu(h_{\alpha}))\]
\end{pf}

\subsection{A formula for $\varphi_{0}$}
We can now derive an explicit formula for the local component 
$\varphi_{0,v}$ of $\varphi_{0}$ and obtain some parity results 
from it. Recall that $v$ has been fixed and we usually don't
include it in the notation. Also, recall the Bruhat decomposition
$\SLT=\widetilde{\B}\amalg\widetilde{\B}w\widetilde{\B}$
\begin{lem} \label{bmatch}
In order for the functions $W_{\varphi}$ and
$W_{\varphi_{0}}$ to agree on $\widetilde{\B}$ it is necessary and
sufficient that
\begin{equation} \label{varphi0}
\varphi_{0}(r)=|r|\chi(r)\int_{T_{x_{0}}(\F_{v})\backslash \HH(\F_{v})}
\varphi(rh^{-1}\cdot x_{0})\chi(\nu(h))dh.
\end{equation}
\end{lem}
\begin{pf}
Since the quadratic form on $\U=\F_{v}$ is $Q(x)=\kappa x^{2}$,
(see Lemma \ref{setm}) 
$W_{\varphi}$ and $W_{\varphi_{0}}$ are in the same Whittaker
space. Therefore, in order for them to agree on $\widetilde{\B}$,
it is necessary and sufficient that
they agree on right $N$ cosets, i.e. on the set $\{ m(a): a\in
\F^{\times} \}$. So we compare
\[
W_{\varphi}(m(a))=|a|^{\frac{3}{2}} (a,-\det(\V))\gamma(a,\psi^{\frac{1}{2}})^{-1}
\int\limits_{\T_{0}\backslash \HH}
\varphi(h^{-1}\cdot a x_{0})\chi(\nu(h))dh
\]
with
\begin{equation*}
W_{\varphi_{0}}(m(a))=|a|^{\frac{1}{2}}(a,2\kappa)\gamma(a,\psi^{\frac{1}{2}})^{-1}\varphi_{0}(a).
\end{equation*}
Solving for $\varphi_{0}(a)$ gives a necessary and sufficient condition for matching.
\[
\varphi_{0}(a)=|a|(a,-\textup{Det}(\V))(a,2\kappa)
\int\limits_{\T_{0}\backslash \HH}
\varphi(h^{-1}\cdot a x_{0})\chi(\nu(h))dh
\]
The quaternion algebra with structure constants $a,b$ has a reduced
norm which can be diagonalized as $(-a,-b,ab)$. This means that the
determinant of the bilinear form associated to the reduced norm is
$8a^{2}b^{2}$ so $(a,-\textup{Det}(\V))(a,2\kappa)=(a,-2)(a,2\kappa)=(a,-\kappa)=\chi(a)$
\end{pf}
In order to to define $\varphi_{0}$ according to the formula in Lemma \ref{bmatch}
we need to show that the formula can be extended to a Schwartz function of $\F_{v}$.
We start with the definition
\begin{defn}
Let $\varphi$ be a Schwartz function in $\Sc(\V)$ and $\chi$ a
quadratic character of $\F_{v}^{\times}$. Define a function on $\F_{v}^{\times}$.
\begin{equation} 
\varphi_{0}^{\times}(r)=|r|\chi(r)\int_{T_{0}\backslash \HH}
\varphi(rh^{-1}\cdot x_{0})\chi(\nu(h))dh
\end{equation}
\end{defn}

\begin{prop}\label{ari}
If $\F_{v}$ is non archimedean, and $\varphi \in \Sc(\V)$ then
$\varphi_{0}^{\times}$ defined by equation \eqref{varphi0} 
can be extended to a function $\varphi_{0}\in \Sc(\U(\F_{v}))$.
\end{prop}

The next few lemmas will establish Proposition \ref{ari}.
The fact that  \ref{varphi0} defines a Schwartz function
at the archimedean place will follow from an indirect
argument in Lemmas \ref{archimatch} and \ref{archimatch2}.

\begin{lem}\label{lc}
When $\F_{v}$ is non archimedean then $\varphi^{\times}_{0}$ is locally
constant 
on $\F_{v}^{\times}$ and supported inside a compact set of $\F$.
\end{lem}
\begin{pf}
Since $\varphi$ has compact support there will be some $n$ such that
$\nu(\textup{supp}(\varphi))\subset \wp^{n}$. Therefore for 
$\rd{r}$ small enough, $\nu(r h^{-1}\cdot x_{0})=r^{2}\kappa \notin
\wp^{n}$ for any $h \in T_{0}\backslash \HH$ and so for $\rd{r}$ small
enough $\varphi_{0}^{\times}(r)=0$. 

Since $\varphi$ is locally constant and has compact support 
there is a neighborhood of $1,\; U\subset \F_{v}^{\times}$ 
such that for all $e\in U$, 
for all $x$ in the support of $\varphi$ we have
$\varphi(ex)=\varphi(x)$. Thus, for $r \neq 0$ and any $e\in U$, 
\begin{align*}
\varphi_{0}(er)&=|er|\chi(r)\int_{\T_{0}\backslash \HH}
\varphi(erh^{-1}\cdot x_{0})\chi(\nu(h))dh\\
&=|r|\chi(r)\int_{\T_{0}\backslash \HH}
\varphi(rh^{-1}\cdot x_{0})\chi(\nu(h))dh\\
&=\varphi_{0}(r)
\end{align*}
\end{pf}

\begin{lem} \label{lcdiv0}
When $\B_{v}$ is the division algebra  ${\mathbb H}$, $\F_{v}$ is non archimedean,
and $\chi$ is nontrivial, then $\varphi_{0}$ vanishes on a neighborhood
of $0$.
\end{lem}
\begin{pf}
Since in the division algebra case $T_{0}\backslash \HH$ 
is compact, for $r$ small enough $\varphi(rh^{-1}\cdot
x_{0})=\varphi(0)$ for all $h$.
For such $r$ 
\[
\varphi_{0}^{\times}(r)=|r|\chi(r)\int_{T_{0}\backslash \HH}
\chi(\nu(h))\varphi(0) dh=0
\]
Since $\ZZ(\F_{v})\backslash \T_{x_{0}}(\F_{v})$ is compact.
and the integral over $T_{0}\backslash \HH$ is equivalent to an
integral over the compact group $\ZZ \backslash \HH$.
\end{pf}

\begin{lem} \label{lcm20}
When $\B_{v}\simeq \M_{2}(\F_{v})$, $\F_{v}$ is non archimedean,
and $\chi$ is non trivial, then $\varphi_{0}$ is constant on a neighborhood
of $0$.
\end{lem}
\begin{pf}
Since $\chi$ is nontrivial, by Lemma \ref{compactness} the set 
$\ZZ(\F_{v})\backslash \T_{x_{0}}(\F_{v})$ is compact and we can define the
function 
\[
\bar{\varphi}(x)=\int_{\ZZ(\F_{v})\backslash \T_{x_{0}}(\F_{v})}
\varphi(t^{-1}x)dt.
\]
It is not hard to see that $\bar{\varphi}$ depends only on the norm of
$x$ and on $\alpha$ as defined in \ref{halpha}.
Next, by Lemma \ref{cosets} we
can see that the measure $dh$ on $\ZZ(\F_{v})\backslash \HH(\F_{v})$ is
a positive scalar multiple of the product measure $d\alpha\, dt$ with
$t \in \ZZ(\F_{v})\backslash \T_{0}(\F_{v})$ and $\alpha$ ranging over
the set of $\alpha$ such that
$\chi(\alpha^{2}-1)=1$. Also, by
Lemma \ref{chihalpha} $\chi(h_{\alpha})=\chi(2(\alpha +1))$.
This means that for the purpose of showing that
$\varphi_{0}$ is constant near $0$ 
we can replace the integral over 
$\ZZ(\F_{v})\backslash \HH(\F_{v})$  by the integral
\[
\varphi_{0}(r)=C\chi(r)|r|
\int_{\{\alpha:\chi(\alpha^{2}-1)=1\}}
\bar{\varphi}(h_{\alpha}^{-1} r x_{0})\chi(2(\alpha +1))d\alpha
\]
for some constant $C$.
With an abuse of
notation we can write $\bar{\varphi}$ as a function of two variables
\[
\bar{\varphi}(\alpha,b)=\int_{\ZZ(\F_{v})\backslash \T_{0}(\F_{v})}
\varphi(t^{-1}(\alpha x_{0}+y))dt.
\]
where $y \perp x_{0}$ and 
$\nu(\alpha x_{0} + y)=b$. Therefore we can write
$\bar{\varphi}(h_{\alpha}^{-1} r x_{0})=\bar{\varphi}(\alpha
r,r^{2}\kappa)$.
Furthermore,
the characteristic function of the set
$\{\alpha:\chi(\alpha^{2}-1)=1\}$ is the function
$\frac{1}{2}[\chi(\alpha^{2}-1)+1]$. When we multiply by the
characteristic function we can extend the integral to an integral over
all of $\F_{v}$. 
\begin{align*}
\varphi_{0}(r)
&=C\chi(r)|r|
\int_{\F_{v}}\bar{\varphi}(\alpha r, r^{2}\kappa)\frac{1}{2}[\chi(\alpha^{2}-1)+1]
\chi(2(\alpha +1))d\alpha \\
&=\frac{\chi(2)C}{2}\chi(r)|r|
\int_{\F_{v}}\bar{\varphi}(\alpha r, r^{2}\kappa)
[\chi(\alpha+1)+\chi(\alpha-1)]d\alpha
\end{align*}
We change variables to $\alpha'=r\alpha$ and we get:
\[
\varphi_{0}(r)=\frac{\chi(2)C}{2}
\int_{\F_{v}}\bar{\varphi}(\alpha', r^{2}\kappa)
[\chi(\alpha'+r)+\chi(\alpha'-r)]d\alpha'
\]
To continue, we will need to extend $x_{0}$ to a basis for $\V$ with 
$x_{0},x_{1},x_{2}$ such that $\nu(x_{0})=\kappa=-\nu(x_{1})$ and
$\nu(x_{2})=-1$ and we also need the following lemma.
\begin{lem}\label{eps}
For some $m>0$,
if $r \in \wp^{m}$ and 
$\chi(\alpha'+r)+\chi(\alpha'-r)\neq 0$ then 
\[
\bar{\varphi}(\alpha', r^{2}\kappa)=
\int_{\ZZ(\F_{v})\backslash \T_{0}(\F_{v})}
\varphi(t^{-1}(\alpha' x_{0}+\alpha' x_{1}))dt
=\bar{\varphi}(\alpha', 0)
\]
\end{lem}
\begin{pf}
The condition that $\chi(\alpha'+r)+\chi(\alpha'-r)\neq 0$ means that 
$\chi(\alpha'^{2}-r^{2})=1$ so the square class of 
$\alpha'^{2}-r^{2}$ can be either $\kappa$ or $1$.

If $\alpha'^{2}-r^{2}=\Box$  i.e. the square class of
$\alpha'^{2}-r^{2}$ is $1$,
let $\epsilon_{1}=-\alpha'
+\sqrt{\alpha'^{2}-r^{2}}$, it is not hard to see that
$|\epsilon_{1}|\leq|\frac{r}{2}|$ regardless of $\alpha'$.
Let $x=\alpha' x_{0}+\alpha' x_{1} + \epsilon_{1}x_{1}$, then
$\nu(x)=r^{2}\kappa $ and the $x_{0}$ coordinate of $x$ is $\alpha'$
so 
\begin{align*}
\bar{\varphi}(\alpha', r^{2}\kappa)
&=\int_{\ZZ(\F_{v})\backslash \T_{0}(\F_{v})}
\varphi(t^{-1}(\alpha' x_{0}+\alpha' x_{1} + \epsilon_{1}x_{1}))dt
\end{align*}
Since $\varphi$ is locally
constant and $\ZZ(\F_{v})\backslash \T_{0}(\F_{v})$ is compact,
for $r$ small enough
\[
\bar{\varphi}(\alpha', r^{2}\kappa)=
\int_{\ZZ(\F_{v})\backslash \T_{0}(\F_{v})}
\varphi(t^{-1}(\alpha' x_{0}+\alpha' x_{1}))dt=\bar{\varphi}(\alpha',0).
\]
Similarly when $\alpha'^{2}-r^{2}=\kappa \Box$ since $\chi$ is not the
trivial character, $\kappa \neq -\Box$ so it is not possible that $|r|>q^{g}|\alpha'|$
where $g$ is the index of the coset of squares in $1 + \wp$, also
since in this case $\alpha'^{2}-r^{2}\neq \Box$, it is not possible that
$|\alpha|>q^{g}|r|$. 

Let $\epsilon_{1}=-\alpha'$ and
$\epsilon_{2}=\sqrt{(\alpha'^{2}-r^{2})\kappa}$. As above,
$|\epsilon_{i}|\leq q^{g} |r|$ and if we let
$x=\alpha' x_{0} + \alpha' x_{1} +\epsilon_{1}x_{1} + \epsilon_{2}
x_{2}$ we get that for small enough $r$
\begin{align*}
\bar{\varphi}(\alpha', r^{2}\kappa)
&=\int_{\ZZ(\F_{v})\backslash \T_{0}(\F_{v})}
\varphi(t^{-1}(\alpha' x_{0}+\alpha' x_{1} +
\epsilon_{1}x_{1}+\epsilon_{2}x_{2}))dt\\
&=\int_{\ZZ(\F_{v})\backslash \T_{0}(\F_{v})}
\varphi(t^{-1}(\alpha' x_{0}+\alpha' x_{1}))dt\\
&=\bar{\varphi}(\alpha',0)
\end{align*}
\end{pf}
To continue with the proof of Lemma \ref{lcm20} notice that for 
$\alpha'$ in some neighborhood of $0$ we have
\[
\int_{\ZZ(\F_{v})\backslash \T_{0}(\F_{v})}
\varphi(t^{-1}(\alpha' x_{0}+\alpha' x_{1}))dt=\varphi(0)
\]
since we normalized the volume of $\ZZ(\F_{v})\backslash \T_{0}(\F_{v})$ to be $1$. 

To show that for $r$ small enough $\varphi_{0}(r)$ is independent of
$r$ let $r$ be small enough that 
$\bar{\varphi}(\alpha', r^{2}\kappa)=\bar{\varphi}(\alpha', 0)$
for any $\alpha'$ such that  $\chi(\alpha'+r)+\chi(\alpha'-r)\neq 0$.
This means that 
\[
\varphi_{0}(r)=\frac{\chi(2)C}{2}
\int_{\F_{v}}\bar{\varphi}(\alpha', 0)
[\chi(\alpha'+r)+\chi(\alpha'-r)]d\alpha'.
\]
If $\alpha'\notin r \wp^{-g}$ then
$\chi(\alpha'-r)=\chi(\alpha'+r)=\chi(\alpha')$
so we can split the integral into two parts
\begin{align*}
\varphi_{0}(r)&=\chi(2)C
\int_{\alpha' \notin r \wp^{-g}}\bar{\varphi}(\alpha', 0)
\chi(\alpha')d\alpha'\\
&+\frac{\chi(2)C}{2}
\int_{r \wp^{-g}}\bar{\varphi}(\alpha', 0)
[\chi(\alpha'+r)+\chi(\alpha'-r)]d\alpha'
\end{align*}
For $r$ small enough the integral over $r\wp^{-g}$ becomes
\[
\frac{\chi(2)C}{2}
\int_{r \wp^{-g}}\varphi(0)
[\chi(\alpha'+r)+\chi(\alpha'-r)]d\alpha'
\]
which splits into the sum of the integrals
\[
\frac{\chi(2)C}{2}\varphi(0)\left(
\int_{r \wp^{-g}}
\chi(\alpha'+r)d\alpha'
\int_{r \wp^{-g}}
+\chi(\alpha'-r)d\alpha' \right)
\]
and when we change variables in the two integrals we get
\[
\chi(2)C\varphi(0)
\int_{r \wp^{-g}}
\chi(\alpha')d\alpha'
\]
so the sum of the two integrals is
\[
\varphi_{0}(r)=\chi(2)C
\int_{\F_{v}}\bar{\varphi}(\alpha', 0)\chi(\alpha')d\alpha'
\]
\end{pf}

We have now established proposition \ref{ari}.

The following proposition is needed in order to
define the global representation.

\begin{prop}[fundamental Lemma] \label{fundi}
If $\F_{v}$ is nonarchimedean with $q$ odd, 
$\chi$ is unramified, $\textup{Inv}(\B_{v})=1$, 
$\nu(x_{0})\in \OO_{v}^{\times}$ 
and $\varphi=\textup{Char}_{\V\cap\M_{2}(\OO_{v})}$ then
$\varphi_{0}=\textup{Char}_{\OO_{v}}$.
\end{prop}
\begin{pf}
From Lemma
\ref{bmatch} we have that when $\chi= 1,\;\varphi_{0}$ must be
\begin{equation} 
\varphi_{0}(r)=|r|\int_{T_{0}(\F_{v})\backslash \HH(\F_{v})}
\varphi(rh^{-1}\cdot x_{0})dh
\end{equation}
If we normalize the measure on $T_{x_{0}}(\F_{v})\backslash
\HH(\F_{v})$ such that $\varphi_{0}(1)=1$ then 
the volume of the set
$\{ h \in T_{0}(\F_{v})\backslash \HH(\F_{v}) : h^{-1}\cdot
x_{0}\in \M_{2}(\OO_{v}) \}$ is one. 

Notice that for $\rd{r}=n$ the set 
$\{ h \in T_{0}(\F_{v})\backslash \HH(\F_{v})\; : \; h^{-1}\cdot
rx_{0}\in \M_{2}(\OO_{v}) \}$ is the same as the set 
$\{ h \in T_{0}(\F_{v})\backslash \HH(\F_{v})\; : \; h^{-1}\cdot
x_{0}\in \varpi^{-n}\M_{2}(\OO_{v}) \}$ 
and so by Lemma \ref{volomegaL} its volume is $q^{n}$. This
means that for $r \in \OO_{v}$ we have $\varphi_{0}(r)=1$. Also
for $r \notin \OO_{v}, \; \nu(h\cdot rx_{0}) \notin \OO_{v}$
so $h\cdot r x_{0} \notin \M_{2}(\OO_{v})$ and 
$\varphi_{0}(r)=0$.

When $\chi$ is unramified non trivial we
normalize the measure such that 
$\textup{Vol}_{\T_{0}\backslash \HH}(\HH(\OO_{v}))=1$
then by Lemma \ref{unramfundi} $\varphi_{0}=\textup{Char}_{\OO_{v}}$.
\end{pf}

\begin{lem} \label{volomegaL}
When $\chi_{v}=1$ (and hence $\B_{v}\simeq \M_{2}$),
$\F_{v}$ is non archimedean and $L \subset \V(\F_{v})$ is
any lattice then if 
\[ 
\textup{Vol}
(\{h \in T_{0}(\F_{v})\backslash \HH(\F_{v}) \; :\;h^{-1} \cdot x_{0}\in L\}) 
 := \textup{Vol}_{\HH}(L) \neq 0
\]
then $\textup{Vol}_{\HH}(\varpi^{-n}L)=q^{n}\textup{Vol}_{\HH}(L)$.
\end{lem}
\begin{pf}
A few things are clear: It is enough to show the result in the case 
of $n=1$, and that
the volume is independent of the choice of $x_{0}$ of norm $\kappa$.
It is also not hard to see that we can assume that $x_{0} \in L$.

To carry out the calculation, pick an isomorphism  of $\B$ with 
$\M_{2}(\F_{v})$ such that $L$ is the trace zero matrices with 
integer coefficients and 
\[ x_{0}=\begin{pmatrix} -1 & 0 \\ 0 & 1 \end{pmatrix} \]

Let $K=\GL_{2}(\OO_{v})$. It turns out that the double coset 
$\T_{0}\backslash \GL_{2}(\F_{v})/K$ is represented by the
set $\{q_{i}: i= 0 \ldots \infty \}$ with
\[
q_{0} = \begin{pmatrix} 1 & 0 \\ 0 & 1 \end{pmatrix}\;\;
q_{i} =  \begin{pmatrix} 1 & \varpi^{-i} \\ 0 & 1 \end{pmatrix}
i > 0.
\]
This means that the coset space
$\T_{x_{0}}\backslash \HH$ is the disjoint union
of the cosets $q_{i}K$.
It is clear that $K \cdot L = L$  and
 $K \cdot \varpi^{n} L = \varpi^{n} L$ so 
$(q_{i}K)^{-1}\cdot x_{0} \subset L$  if and only if $q_{i}^{-1}\cdot
x_{0}\in L$ and also $(q_{i}K)^{-1}\cdot x_{0} \subset \varpi^{-1} L$
 if and only if $q_{i}^{-1}\cdot x_{0}\in \varpi^{-1} L$.

The only $q_{i}$ such that  
$q_{i}\cdot x_{0} \in L$ is $q_{0}$, so that 
$\textup{Vol}_{\HH}(L)$ is equal to the volume of the set 
$K \subset \T_{0}\backslash \HH$. 

Similarly $q_{1}$ and $q_{0}$ are the only  $q_{i}$ such that  
$q_{i}\cdot x_{0} \in \varpi^{-1} L$ and so 
$\textup{Vol}_{\HH}(L)$
 is the sum of the volumes of 
$K$ and $q_{1}K$. Comparison of the volumes of the two lattices now 
follows from a simple group calculation.
\end{pf}

\begin{lem} \label{unramfundi}
If $q$ is odd, $\chi$ is unramified non trivial,
 $\varphi={\textup Char}(\M_{2}(\OO_{v}))$
and $\textup{Vol}_{\T_{0}\backslash \HH}(\M_{2}(\OO_{v}))=1$ then
$\varphi_{0}$ as defined by equation \ref{varphi0} 
is equal to $\textup{Char}(\OO_{v})$.
\end{lem}
\begin{pf}
Pick 
\[
x_{0}=\begin{pmatrix}0 & -\gamma \\ 1 & 0 \end{pmatrix}
\]
with $\gamma \neq -\Box,\;\;\gamma \in \OO_{v}^{\times}$. By Lemma
\ref{ext} we can see that
\[
\T_{0}=\left\{ \begin{pmatrix}a & -\gamma b \\ b & a \end{pmatrix} \;:\;a,b \in \F_{v}  \right\}
\]
Notice that the function $\varphi$ is invariant by the action of 
$K=\HH^{\times}(\OO_{v})=\GL_{2}(\OO_{v})$ so to calculate
\[
\varphi_{0}(r)=\chi(r)|r|\int_{\T_{0}\backslash
\HH}\varphi(h^{-1}\cdot r x_{0})\chi(\nu(h))dh
\]
The integral can be replaced by a sum over coset representative
of 
\[
\T_{0}\backslash \HH/ K = 
\left\{ \begin{pmatrix} 1 & 0 \\ 0 & \varpi^{n}\end{pmatrix}:=q_{n}\; :\; n\geq 0 \right\}
\]
To actually calculate the integral we need to calculate the 
volume of the cosets $\textup{Vol}_{\T_{0}\backslash \HH}(q_{i}K)$. 
Since we set $\textup{Vol}_{\T_{0}\backslash \HH}(K)=1$ we obtain 
the volume of the other cosets by comparing them to the
volume of $K$. And in fact, for $i>0$,
$\textup{Vol}_{\T_{0}\backslash \HH}(q_{i}K)=(q+1)q^{i-1}$.
We can now calculate $\varphi_{0}$. First notice that if 
$r\notin \OO_{v}$ then $\nu(rx_{0})\notin \OO_{v}$ and $\varphi_{0}(r)=0$.
When $r\in \OO_{v}$ then 
\[
\varphi_{0}(r)=
|r|\chi(r)\left( \varphi(r x_{0})
+\sum_{i=1}^{\infty}(q+1)q^{i-1}
\varphi(q_{i}^{-1}\cdot x_{0})\chi(\nu(q_{i})) \right)
\]
now, $q_{i}^{-1}\cdot r x_{0} =
\begin{pmatrix} 0 & -\varpi^{-i}r\gamma \\ r\varpi^{i} & 0 \end{pmatrix}$
so for $i<\rd{r},\; q_{i}^{-1}\cdot r x_{0}$ is not in the support of
$\varphi$ so
\[
\varphi_{0}(r)
=q^{-\rd{r}}(-1)^{\rd{r}}\left( 1
+\sum_{i=1}^{\rd{r}}(q+1)q^{i-1}(-1)^{i} \right)
=1
\]
\end{pf}
We can now obtain some parity properties for $\varphi$ and $\varphi_{0}$.
\begin{cor}\label{par}
The parity of $\varphi_{0}$ is even when $\textup{Inv}(\B_{v})=1$
and odd if $\textup{Inv}(\B_{v})=-1$.
\end{cor}

\begin{pf}
Be Lemma \ref{conji} there is an element $\xi \in \HH(\F_{v})$ 
such that $\xi^{-1}\cdot x_{0}=-x_{0}$ and $\chi(-\nu(\xi))=\textup{Inv}(\B_{v})$.
With that we can calculate 
\begin{align*}
\varphi_{0}(-r)&=|r|\chi(-r)\int_{T_{x_{0}}\backslash \HH}
\varphi(r(\xi h)^{-1}) \cdot x_{0})\chi(\nu(h))dh \\
&=|r|\chi(-r)\int_{T_{x_{0}}\backslash \HH}
\varphi(r h'^{-1} \cdot x_{0})\chi(\nu(\xi^{-1}h'))dh' 
\end{align*}
with $h'=\xi h$
\begin{align*}
&=|r|\chi(-r)\chi(\nu(\xi))\int_{T_{x_{0}}\backslash \HH}
\varphi(r(h'^{-1}) \cdot x_{0})\chi(\nu(h'))dh'\\
&=\chi(-\nu(\xi))|r|\chi(r)\int_{T_{x_{0}}\backslash \HH}
\varphi(r(h'^{-1})  \cdot x_{0})\chi(\nu(h'))dh'\\
&=\textup{Inv}(\B)\varphi_{0}(r)
\end{align*}
\end{pf}

\begin{cor}\label{par2}
Let $\varphi \in \Sc(\V(\F_{v}))\; \varphi=\varphi^{+}+\varphi^{-}$
with $\varphi^{+}$ an even function and $\varphi^{-}$ an odd function.
If $\chi(-1)\textup{Inv}(\B_{v})=1$ then $\varphi^{-}_{0}=0$ and 
if $\chi(-1)\textup{Inv}(\B_{v})=-1$ then $\varphi^{+}_{0}=0$.
\end{cor}
\begin{pf}Consider $\varphi^{+}$.
By corollary \ref{par}
$\varphi_{0}(-r)=\textup{Inv}(\B_{v})\varphi_{0}(r)$ so
\begin{align*}
\varphi_{0}^{+}(r)
&=\textup{Inv}(\B_{v})\varphi_{0}^{+}(-r)\\
&=\textup{Inv}(\B_{v})|r|\chi(-r)\int_{T_{0}\backslash \HH}
\varphi^{+}(-rh^{-1}\cdot x_{0})\chi(\nu(h))dh\\
&=\textup{Inv}(\B_{v})\chi(-1)|r|\chi(r)\int_{T_{0}\backslash \HH}
\varphi^{+}(rh^{-1}\cdot x_{0})\chi(\nu(h))dh\\
&= \textup{Inv}(\B_{v})\chi(-1)\varphi_{0}^{+}(r)
\end{align*}
and we see that if $\textup{Inv}(\B_{v})\chi(-1)\neq 1$
then $\varphi_{0}^{+}(r)=-\varphi_{0}^{+}(r)$. The case of $\varphi_{0}^{-}$ is similar.
\end{pf}

Recall that the choice of parity at places of $\F$ determines the
character $\mu$. We can restate corollary \ref{par} in terms of $\mu$
and obtain the remaining part of the local data defining $W(\U,\mu)$.

\begin{lem}
In order that $W_{\V}(\chi)\subset W_{\U}(\mu)$
it is necessary that the character $\mu=\otimes_{v}\mu_{v}$ is
defined by
\[
\mu_{v}=
\begin{cases}
1 &\textup{ if } \textup{Inv}(\B_{v})=1 \\
\textup{Sgn}& \textup{ if }   \textup{Inv}(\B_{v})=-1
\end{cases}
\]
\end{lem}

\section{Nonarchimedean local matching}

\subsection{An equivariant property}
For the following, we will need to recall
the Bruhat decomposition.
From Lemma \ref{bmatch} we know that if $\varphi_{0}$ is defined according to
equation \eqref{varphi0} then $W_{\varphi}$ and $W_{\varphi_{0}}$ agree
on $\widetilde{B}$. It remains to show that with the same
$\varphi_{0}$, the Whittaker functions $W_{\varphi}$ and
$W_{\varphi_{0}}$ will also agree on $\widetilde{B}w\widetilde{B}$.
The verification of this result
 is simplified by the following lemma.
\begin{lem} \label{coin}
Let the map $i: \Sc(\V)\rightarrow \Sc(\U)$ be defined by $i(\varphi)=\varphi_{0}$
 where
$\varphi_{0}$ is given by 
equation \eqref{varphi0}, then for any $g \in \tilde{B}$
\[
i(\omega(g)(\varphi))=\omega_{0}(g)(i(\varphi)).
\]
\end{lem}
\begin{pf}
Since $ \tilde{B}=NA=AN$, it will suffice to
show that for any $a\in \F_{v}^{\times}$, 
\begin{equation}\label{aaa}
i(\omega(m(a))(\varphi))=\omega_{0}(m(a))(i(\varphi))
\end{equation}
and for any $b\in \F_{v}$
\begin{equation}\label{bbb}
i(\omega(n(b))(\varphi))=\omega_{0}(n(b))(i(\varphi)).
\end{equation}
To establish equation $\eqref{aaa}$ we compute $i(\omega(m(a))(\varphi))(r)$
\begin{align*}
&=|r|\chi(r)\int_{T_{0}\backslash \HH}
\omega(m(a))\varphi(rh^{-1}\cdot x_{0})\chi(\nu(h))dh\\
&=|r|\chi(r)|a|^{\frac{3}{2}}(a,-\det(\V))\gamma(a,\psi^{\frac{1}{2}})^{-1} 
\int_{T_{0}\backslash \HH}
\varphi(arh^{-1}\cdot x_{0})\chi(\nu(h))dh.
\end{align*}
On the other hand $\omega_{0}(m(a))(i(\varphi))(r)$
\begin{align*}
&=|a|^{\frac{1}{2}}(a,2\kappa)\gamma(a,\psi^{\frac{1}{2}})^{-1}|ar|\chi(ar)
\int_{T_{x_{0}}\backslash \HH}
\varphi(arh^{-1}\cdot x_{0})\chi(\nu(h))dh\\
&=|r|\chi(r)|a|^{\frac{3}{2}}\chi(a)(a,2\kappa)\gamma(a,\psi^{\frac{1}{2}})^{-1}
\int_{T_{0}\backslash \HH}
\varphi(arh^{-1}\cdot x_{0})\chi(\nu(h))dh.
\end{align*}
Comparing the two sides of equation \eqref{aaa}, we see that they are equal if
$\chi(a)(a,2\kappa)=(a,-\det(\V)).$
Recall from proof of lemma \ref{bmatch}
 that $\det(\V))$ is in the square class
of $2$ for both division and matrix algebra so it is true that
$\chi(a)(a,2\kappa)=(a,-\kappa)(a,2\kappa)=(a,-\det(\V)).$

To establish equation \eqref{bbb} we calculate 
 $i(\omega(n(b))(\varphi))(r)$
\begin{align*}
&=|r|\chi(r)\int_{T_{0}\backslash \HH}
\omega(n(b))\varphi(rh^{-1}\cdot x_{0})\chi(\nu(h))dh\\
&=\psi(b\kappa)|r|\chi(r)
\int_{T_{0}\backslash \HH}
\varphi(rh^{-1}\cdot x_{0})\chi(\nu(h))dh\\
&=\psi(b\kappa)i(\varphi)(r)\\
&=\omega_{0}(n(b))i(\varphi))(r)
\end{align*}
\end{pf}

\begin{lem} \label{coinvariant}
If  for all $ \varphi \in \Sc(\V),\;\;
W_{\varphi}(w)=CW_{\varphi_{0}}(w)$ for some constant $C$, then  
for all $\varphi \in \Sc(\V)$ and for all $g\in \widetilde{B}w\widetilde{B}$,
$W_{\varphi}(g)=CW_{\varphi_{0}}(g)$.
\end{lem}
\begin{pf}
The subgroup $\tilde{B}$ can be factored as $NA$.
Since $m(a)w=wm(a^{-1})$
The cell $\widetilde{B}w\widetilde{B}$ can be written as
$Nw\widetilde{B}$.
Consider any $wg \in wB$, by assumption, for any $\varphi\in \Sc(\V))$
we have $W_{\varphi}(w)=CW_{i(\varphi)}(w)$ so, using Lemma \ref{coin}
\begin{align*}
W_{\varphi}(wg)
&=W_{\omega(g)\varphi}(w)\\
&=CW_{i(\omega(g)\varphi)}(w)\\
&=CW_{\omega_{0}(g)i(\varphi)}(w)\\
&=CW_{i(\varphi)}(wg)
\end{align*}
So $W_{\varphi}$ and $CW_{i(\varphi)}$ coincide on all of $w\tilde{B}$.
Finally, since both $W_{\varphi}(w)$ and $CW_{i(\varphi)}(w)$ are in
the same Whittaker space by construction, they coincide on all of
$Nw\tilde{B}=\tilde{B}W\tilde{B}.$
\end{pf}

We will use two different arguments to show that $W_{\varphi}$ matches
$W_{\varphi_{0}}$. In some cases we will match a single $\varphi$ and
use general information about the Weil representation to conclude that 
in fact, all other functions $\varphi$ are also matched by
$\varphi_{0}$.
In other cases we will show directly that for any $\varphi \in \Sc(\V),\;\;
W_{\varphi}(w)=CW_{\varphi_{0}}(w)$ for some constant $C$ and then use Lemma
\ref{coinvariant} to conclude that the
$W_{\varphi}=CW_{\varphi_{0}}$ on all of $\tilde{B}w\tilde{B}$ for
some constant $C$. Finally, we will show that $C=1$ and 
deduce the actual equality  on all of $\SLT$
from that relation.

\subsection{$\chi\neq 1$.}
 From equation \eqref{w0w} and the definition 
of $W_{\varphi_{0}}$, ( definition
\ref{wvarphi0}.) we have that
\[
W_{\varphi_{0}}(w)=\gamma(\psi\circ \U)^{-1}\int_{\U} \varphi_{0}(y)\psi(-2\kappa y)dy
\]
Notice that we set $m$ in \eqref{w0w} equal to $\kappa$ according to Lemma \ref{setm}.

The measure $dy$ is self dual with respect to the Fourier transform
$\hat{f}(x)=\int_{\U(\F_{v})}f(y)\psi(-2\kappa x y)dy$. Since
the conductor of $\psi$ is $\epsilon^{-1}\varpi^{-d}\OO_{v}$ it follows by a
standard argument that $dy=|\varpi^{d}\epsilon 2\kappa|^{\frac{1}{2}}du$ where
$du$ gives $\OO_{v}$ volume $1$.

Since $\det(\U)=2\kappa$ and $\epsilon(\U)=1$ we can simplify the
expression for
$\gamma(\psi\circ \U)=\gamma(\det(\U),\psi^{\frac{1}{2}})\gamma(\psi^{\frac{1}{2}})\epsilon(\U)$
to $\gamma(\psi^{\kappa})$. And substituting the explicit expression
for $\varphi_{0}$ in terms of $\varphi$, we have
\begin{align} 
W_{\varphi_{0}}(w)=&|\varpi^{d}\epsilon 2 \kappa|^{\frac{1}{2}}
\gamma(\psi^{\kappa})^{-1}\int_{\U}
\varphi_{0}(u)\psi(-2\kappa u)du \nonumber \\
=&|\varpi^{d}\epsilon 2\kappa|^{\frac{1}{2}}
\gamma(\psi^{\kappa})^{-1}\int_{\U}|u|\chi(u)
\nonumber \int_{T_{0}\backslash \HH}\!\!\!\!\!
\varphi(h^{-1}\cdot ux_{0})\nonumber \\
&\hspace{85pt}\times \chi(\nu(h))dh\;\psi(-2\kappa u)du\label{wvarphi0w}
\end{align}

On the other hand, by the definition of the action of $\omega(w)$
(equation \eqref{ww}) and the definition of $W_{\varphi}$  we have
\begin{equation}\label{wvarphiw}
W_{\varphi}(w)= \gamma(\psi\circ \V)^{-1}\int\limits_{\T_{0}\backslash \HH}
\int_{\V} \varphi(y)\psi(-(h^{-1}x_{0},y))dy\,\chi(\nu(h))dh
\end{equation}

\begin{rem}
As mentioned above, we need to verify the equality of
\eqref{wvarphi0w} and \eqref{wvarphiw}. Notice that both equations 
involve integrals over $\V(\F_{v})$ and $T_{0}\backslash \HH$
so verification of the equality is independent of 
the realization of the local algebra $\B_{v}$
as long as the measures $dy$ and $du$ are self dual with respect to
the action of $w$ and
the measure on $T_{0}\backslash \HH$ is the same in
both equations.
\end{rem}

We continue with the calculation of the right hand side of
\eqref{wvarphiw}. Since $\chi(r)=(r,-\kappa)$ is nontrivial, and $\kappa=\nu(x_{0})$ we know from
Lemma \ref{compactness} that $\ZZ \backslash \T_{0}$ is compact.
It is also clear that the integrand in \eqref{wvarphiw} is constant on left $\T_{0}$-cosets
so we will replace the integral over 
$\T_{0}\backslash \HH $ with an integral over
$Z\backslash \HH$. To do that,
Let $M_{\kappa}=\textup{Vol}(\ZZ \backslash \T_{0})$ and we have:
\begin{align*}
W_{\varphi}(w)
&= \gamma(\psi\circ \V)^{-1}\int\limits_{\T_{0}\backslash \HH}
\int_{\V} \varphi(y)\psi(-(h^{-1}x_{0},y))dy\,\chi(\nu(h))dh \\
&= M_{\kappa}^{-1}\gamma(\psi\circ \V)^{-1}\int\limits_{\ZZ \backslash \HH}
\int_{\V} \varphi(y)\psi(-(h^{-1}x_{0},y))dy\,\chi(\nu(h))dh.
\end{align*}
We now change variables in the inner integral from $y$ to $y'=h\cdot y$
\begin{align*}
W_{\varphi}(w)
&= \gamma(\psi\circ \V)^{-1}M_{\kappa}^{-1}\int\limits_{\ZZ \backslash \HH}
\int_{\V}\varphi(h^{-1}y')\psi(-(x_{0},y'))dy'\,\chi(\nu(h))dh. 
\end{align*}
Next we let $C_{L}$ be the characteristic function of a compact set
$L\subset \V(\F_{v})$
and we note that for any increasing sequence of compact sets  
\begin{align*}
&W_{\varphi}(w)\\
&= \lim_{L}\gamma(\psi\circ \V)^{-1}M_{\kappa}^{-1}\!\int\limits_{\ZZ \backslash \HH}
\!\int_{\V} \varphi(h^{-1}y)\psi(-(x_{0},y))C_{L}(y)dy\,\chi(\nu(h))dh \\
&= \lim_{L}\gamma(\psi\circ \V)^{-1}M_{\kappa}^{-1}\!\int_{\V}\psi(-(x_{0},y))
\!\int\limits_{\ZZ \backslash \HH}
\varphi(h^{-1}y)\chi(\nu(h))dh\,C_{L}(y)dy\\
&=\lim_{L}\gamma(\psi\circ \V)^{-1}M_{\kappa}^{-1}\!\int_{\om{}^{\times}}\psi(-(x_{0},y))
\!\int\limits_{\ZZ \backslash \HH}
\varphi(h^{-1}y)\chi(\nu(h))dh\,C_{L}(y)dy,
\end{align*}
where in the last step
we replaced the integral over $\V(\F_{v})$ by the
integral over $\om{}^{\times}=\{v\in \V \,:\, \nu(v)\neq 0 \}$.
This is justified by the fact that the set of elements with $0$ norm
have measure $0$ in $\V$.

Now, from Lemma \ref{ker} we know that if $\nu(y)$ is not in the
square class of $\kappa$ then there is an element $h_{y}$ in its stabilizer
$T_{y}(\F_{v})$ such that $\chi(\nu(h_{y}))=-1$. In that case the
inner integral is zero. So we can assume that the support of $\varphi$
is contained in the set 
\[
\om{\kappa}^{\times}=\{v\in \om{}^{\times} \,:\, \nu(v)\in 
\kappa(\F_{v}^{\times})^{2}\}
\]
We state this as a lemma for reference.
\begin{lem} \label{sup}
Let $\chi(a)=(a,-\kappa)$ be a nontrivial character and let
$\varphi_{\kappa}$ be the restriction of $\varphi$ to $\om{\kappa}^{\times}$
then $W_{\varphi}(w)=W_{\varphi_{\kappa}}(w)$.
\end{lem}

To continue, we will parameterize $\om{\kappa}^{\times}$ by using the two fold
covering 
\[
\F^{\times} \times \T_{0}\backslash \HH \rightarrow
\om{\kappa}^{\times} \;\;\;\;\; (c,h) \rightarrow ch^{-1}x_{0}.
\]

The measure on $\om{}^{\times}$ is the measure $dy$ which makes the
Fourier transform $\hat{f}(x)=\int_{\V(\F_{v})}\psi(-(x,y))f(y)dy$
self dual. This measure is invariant under the action of $\HH$
and transforms by $|c|^{3}$ under multiplication by a scalar $c$.

The measure $|u|^{2}dhdu$ on the cover 
$\F^{\times}\times \T_{x_{0}}\backslash \HH$ is also invariant under 
the action of $\HH$ and also transforms by $|c|^{3}$ under multiplication
by $c$.
It follows that the two measures are related by a positive constant
$C$. We defer the calculation of $C$ since we will soon change
variables again and will normalize the measure that results from the
two changes of variables.

The relation between the two measures allows us to write
$W_{\varphi}(w)$ in terms of integrals over 
$\F^{\times}\times \T_{0}\backslash \HH$ 
\begin{align*}
W_{\varphi}(w)=&\lim_{L}\gamma(\psi \circ V)^{-1}CM_{\kappa}^{-1}
\int\limits_{\F_{v}^{\times}}
\int\limits_{\T_{0}\backslash \HH}
\psi(-(x_{0},h_{1}^{-1}ux_{0}))\\
&\times \int_{\ZZ\backslash \HH}
\varphi(h^{-1}h_{1}^{-1}ux_{0})\chi(\nu(h))dh\,C_{L}(h_{1}^{-1}ux_{0})
|u|^{2}dh_{1}\,du 
\end{align*}
We let $h_{2}=h_{1}h$ and we get
\begin{align*}
W_{\varphi}(w)=\lim_{L}\gamma(\psi \circ &V)^{-1} CM_{\kappa}^{-1}
\int\limits_{\F_{v}^{\times}}
\int\limits_{\T_{0}\backslash \HH}
\psi(-(x_{0},h_{1}^{-1}ux_{0}))\chi(\nu(h_{1}))|u|^{2}\\
&\times \int_{\ZZ\backslash \HH}
\varphi(h_{2}^{-1}ux_{0})\chi(\nu(h_{2}))dh_{2}\;C_{L}(h_{1}^{-1}ux_{0})
dh_{1}\,du 
\end{align*}

Notice that the integrand in the integral over
$\ZZ\backslash \HH$ is left invariant under 
$\ZZ\backslash \T_{0}$ 
and recall that the volume of 
$\ZZ\backslash \T_{0}$ is $M_{\kappa}$
and so the inner integral can be changed
to an integral over
$\T_{0}\backslash \HH$.
\begin{align*}
W_{\varphi}(w)=&\lim_{L}\gamma(\psi \circ V)^{-1} C
\int\limits_{\F_{v}^{\times}}\\
&
\times\left(
\int\limits_{\T_{0}\backslash \HH}
\psi(-(x_{0},h_{1}^{-1}ux_{0}))\chi(\nu(h_{1}))|u|^{2} C_{L}(h_{1}^{-1}ux_{0})dh_{1}\right)    \\
&\times \left( \int_{\T_{0}\backslash \HH}
\varphi(h_{2}^{-1}ux_{0})\chi(\nu(h_{2}))dh_{2}\right)du.
\end{align*}

Comparing with equation \eqref{varphi0} we see that 
the third integral is $\varphi_{0}(u)|u|^{-1}\chi(u)$ so we get
\begin{align*}
W_{\varphi}(w)=&\lim_{L}\gamma(\psi \circ V)^{-1} C
\int\limits_{\F_{v}^{\times}}\\
&\times\left( \int\limits_{\T_{0}\backslash \HH}
\psi(-(x_{0},h_{1}^{-1}ux_{0}))\chi(\nu(h_{1})) 
 C_{L}(h_{1}^{-1}ux_{0})dh_{1} \right)\\
&\times |u|\chi(u)\varphi_{0}(u) du.
\end{align*}

The integrand in $W_{\varphi}(w)$ depends only on cosets of the form 
$h_{\alpha}T_{0}$ so using Lemmas
 \ref{cosets}
\ref{halpha} we can change the integral in
$W_{\varphi}(w)$ over $\T_{0}\backslash \HH $ to
an integral over the product of the set of $\alpha$ such that 
$\chi(\alpha^{2}-1)=\textup{Inv}(\B_{v})$ and $\ZZ(\F_{v})\backslash
\T_{x_{0}}(\F_{v})$. 

We have now changed measures from $dy$ to the product of the
measures $d\alpha\, du\, dt$, with $\alpha$ as above, $t \in
\ZZ\backslash \T_{0}$ and $u \in \F_{v}^{\times}$. At this point in the
calculation we note that  
the constant relating the two measures
will depend on the type of 
algebra $\B_{v}$ and 
on the norm of $x_{0}$ defining the torus $\T_{0}$.
We denote the constant relating the measures $dy$ and  $|u^{2}|d\alpha\, du\, dt$
by $C_{\kappa,v}$.

So far
the formula for $W_{\varphi}(w)$ looks like:

\begin{align} 
C_{\kappa,v}\lim_{L}\gamma(\psi \circ V)^{-1}&
\!\!\int\limits_{\F_{v}^{\times}}\!|u|\chi(u) \varphi_{0}(u)
\!\!\int\limits_{\alpha}\!
\int\limits_{\ZZ\backslash \T_{0}}
\!\!\!\!\!
\psi(-(x_{0},t^{-1}h_{\alpha}^{-1}ux_{0}))\chi(\nu(th_{\alpha}))\nonumber \\
&
C_{L}(t^{-1}h_{\alpha}^{-1}ux_{0})dt\,d\alpha\,du, \label{lookslike}
\end{align}
where $\alpha$ ranges over elements such that $\chi(\alpha^{2}-1)=\textup{Inv}(\B)$. 
Also
recall that $h_{\alpha}^{-1}x_{0}=\alpha
x_{0}+ x$ with $x \perp x_{0}$ so
\[ (x_{0},h_{\alpha}ux_{0})=\alpha u
(x_{0},x_{0})= \alpha u 2 \nu(x_{0})= 2\alpha u \kappa
\]
With the above observation and with Lemmas \ref{cosets} and \ref{chihalpha}
the two inner integrals become
\[C_{\kappa,v}\lim_{L}\int\limits_{\alpha}
\int\limits_{\ZZ\backslash \T_{0}}
\!\!\!
\psi(-2u\alpha \kappa )\chi(2(\alpha+1)))
C_{L}(t^{-1}h_{\alpha}^{-1}ux_{0}) dt\,d\alpha .
\]

Since we are taking the limit over all increasing sequences of
lattices $L$ the limit would not change if we considered only lattices
which are invariant under the action of the compact group
$\ZZ\backslash \T_{0}$. In that case, the integrand does
not depend on $t$. Let 
$V_{x_{0}}=\textup{Vol}(\ZZ\backslash \T_{0})$ 
and the inner integrals become
\[
V_{x_{0}}\lim_{L}\int\limits_{\alpha}\psi(-2u\alpha \kappa )
\chi(2(\alpha+1))C_{L}(h_{\alpha}^{-1}ux_{0})
d\alpha .
\]

The characteristic function of the set of $\alpha \neq \pm 1$ such
that  $\chi(\alpha^{2}-1)=\textup{Inv}(\B)$ is the function
\[
\frac{1}{2}(1 + \textup{Inv}(\B)\chi(\alpha^{2}-1))
\]
so we can extend the integral over $\alpha$ to an integral over 
$\F_{v} - \{\pm 1\}$ by multiplying the integrand by the
characteristic function. The integral over $\F_{v}$ is only defined as
a limit and should be understood that when we write
$\int_{\F_{v}}\psi(-\alpha')\chi(\alpha')d\alpha$ we will mean 
$\lim_{n}\int_{\wp^{-n}}\psi(-\alpha')\chi(\alpha')d\alpha$.

\begin{align*}
&V_{x_{0}}\int\limits_{\F_{v}-\{\pm 1\}}\psi(-2u\alpha \kappa )
\chi(2(\alpha+1))\frac{1}{2}(1 + \textup{Inv}(\B)\chi(\alpha^{2}-1))d\alpha \\
&=V_{x_{0}}\chi(2)\int\limits_{\F_{v}-\{\pm 1\}}\psi(-2u\alpha \kappa )
\frac{1}{2}(\chi(\alpha +1) + \textup{Inv}(\B)\chi(\alpha-1)) d \alpha .
\end{align*}
We split the integral into two integrals and change variables to
$\alpha +1$ and to $\alpha -1$ respectively. We also define
$\chi(0)=0$ and extend the integral to all of $\F_{v}$.
\[=V_{x_{0}}\chi(2)\frac{1}{2}(\psi(2\kappa u) + \textup{Inv}(\B)\psi(-2\kappa u))
\int\limits_{\F_{v}}\psi(-2u\alpha \kappa )\chi(\alpha) d\alpha 
\]
We use another change of variable $\alpha' = 2u\alpha \kappa$ and we get
\[
=V_{x_{0}}\chi(u)|2u\kappa|^{-1}\frac{1}{2}(\psi(2\kappa u) + \textup{Inv}(\B)\psi(-2\kappa u))
\int\limits_{\F_{v}}\psi(-\alpha' )\chi(\alpha')
d\alpha' .
\]
We give the above integral a name 
\begin{defn} \label{gauss}
The 'gauss integral' is the quantity
\[
\mathfrak{g}_{v}(\psi,\chi)=\int\limits_{\F_{v}}\psi(-\alpha)\chi(\alpha)d\alpha
=\lim_{n}\int\limits_{\wp^{-n}}\psi(-\alpha)\chi(\alpha)d\alpha
\]
\end{defn} 
So far the inner integrals in equation \ref{lookslike}
are
\begin{equation} \label{sofar1}
\mathfrak{g}_{v}(\psi,\chi)V_{x_{0}}
\chi(u)|2u\kappa|^{-1}\frac{1}{2}(\psi(2\kappa u) +
\textup{Inv}(\B)\psi(-2\kappa u))
\end{equation}

putting equation \ref{sofar1} together with equation \eqref{lookslike}
we get $W_{\varphi}(w)=$
\begin{equation*} 
C_{\kappa,v}\gamma(\psi \circ V)^{-1}|2\kappa|^{-1}
\mathfrak{g}_{v}(\psi,\chi)V_{x_{0}}
\int\limits_{\F_{v}^{\times}} \varphi_{0}(u)
\frac{1}{2}(\psi(2\kappa u) +
\textup{Inv}(\B)\psi(-2\kappa u)) du 
\end{equation*}

Now by Lemma \ref{par}
$\varphi_{0}(-u)=\textup{Inv}(\B_{v})\varphi_{0}(u)$  
and so 
\begin{align*} 
&W_{\varphi}(w)\\
&=C_{\kappa,v}\textup{Inv}(\B_{v})
\gamma(\psi \circ V)^{-1}|2\kappa|^{-1}\mathfrak{g}_{v}(\psi,\chi)V_{x_{0}}
\int\limits_{\F_{v}^{\times}} \varphi_{0}(u)\psi(-2\kappa u) du \\
&=C_{\kappa,v}V_{x_{0}}\mathfrak{g}_{v}(\psi,\chi)\textup{Inv}(\B_{v})
\gamma(\psi \circ V)^{-1}|2\kappa|^{-1}
|\varpi^{d}\epsilon 2\kappa|^{-\frac{1}{2}}
\gamma(\psi^{\kappa})W_{\varphi_{0}}(w)\\
&=C_{\kappa,v}
V_{x_{0}}\mathfrak{g}_{v}(\psi,\chi)
|\varpi^{d}\epsilon|^{-\frac{1}{2}}|2\kappa|^{-\frac{3}{2}}
\textup{Inv}(\B_{v})
\gamma(\psi \circ V)^{-1}\gamma(\psi^{\kappa})W_{\varphi_{0}}(w)
\end{align*}
What we showed so far is that on the ''big cell``, 
the two Whittaker functions $W_{\varphi}$ and $W_{\varphi_{0}}$ 
are related by a constant which
depends on $\chi,\F_{v},\B_{v}$ but does not depend on the argument 
of the function. I.e.
\[
W_{\varphi}(g)=C_{\kappa,v,\B_{v}}W_{\varphi_{0}}(g)
\]
with
\[
C_{\kappa,v,\B_{v}}=
C_{\kappa,v}
V_{x_{0}}\mathfrak{g}_{v}(\psi,\chi)
|\varpi^{d}\epsilon|^{-\frac{1}{2}}|2\kappa|^{-\frac{3}{2}}
\textup{Inv}(\B_{v})
\gamma(\psi \circ V)^{-1}\gamma(\psi^{\kappa})
\]
If we showed that $C_{\kappa,v,\B_{v}}=1$ we would be done.

\begin{prop} \label{nantm}
When $\F_{v}$ is non archimedean and $\chi\neq 1$ then for any $\varphi
\in \Sc(\V)$
\[
W_{\varphi}(w)=W_{\varphi_{0}}(w).
\]
Furthermore, $\Theta(\HH,\chi)\simeq \Theta(\GO(\U),\mu)$ and for any
$\varphi \in \Sc(\V)$ we have $W_{\varphi}=W_{\varphi_{0}}$.
\end{prop}
\begin{pf}
The discussion above
shows that for any $\varphi$, and for $\varphi_{0}$ defined by
equation \eqref{varphi0} $W_{\varphi}(w)=C_{\kappa,v,\B_{v}}W_{\varphi_{0}}(w)$.
By Lemma \ref{coinvariant} we have that 
 $W_{\varphi}(g)=C_{\kappa,v,\B_{v}}W_{\varphi_{0}}(g)$ for all
$g \in \tilde{B}w\tilde{B}$. Since the Weil representation of $\SLT$
is smooth, the stabilizers of $\varphi$ and $\varphi_{0}$ are both
open subgroups and so the intersection of the stabilizers $G$ is an
open neighborhood of $I\in \SLT$. It is not hard to see that
\[
\forall g \in G,\;\;\;W_{\varphi}(g)=W_{\omega(g)\varphi}(I)=W_{\varphi}(I)
\] and
\[
\forall g \in G,\;\;\;W_{\varphi_{0}}(g)=W_{\omega_{0}(g)\varphi_{0}}(I)=W_{\varphi_{0}}(I)
\]
We can find $g'\in G$ such that
$g'\in \tilde{B}w\tilde{B}$. If we pick $\varphi$ such that
$W_{\varphi}(I)\neq 0$ then
\begin{align*}
W_{\varphi}(g')
&=C_{\kappa,v,\B_{v}}W_{\varphi_{0}}(g')\\
&=C_{\kappa,v,\B_{v}}W_{\varphi_{0}}(I)
\end{align*}
and 
\begin{align*}
W_{\varphi}(g')&=W_{\varphi}(I)\\
&=W_{\varphi_{0}}(I)
\end{align*}
since $I\in\tilde{B}$. Therefore $C_{\kappa,v,\B_{v}}=1$ and we are done.
\end{pf}

\subsection{$\chi=1$.}
The argument for matching in this case is different from the argument 
in the case of a nontrivial character. In this
case, the equality of the spaces $W_{\psi^{\kappa}}(\V,\chi)$ and 
$W_{\psi^{\kappa}}(\U,\mu)$ follows from results of Waldspurger and
Rallis. The explicit matching of the functions $W_{\varphi}$ and 
$W_{\varphi_{0}}$ will then follow from Lemma \ref{bmatch}. 

\begin{prop} \label{natm}
For $\chi=1$ and $\F_{v}$ non archimedean
\[
W_{\psi^{\kappa}}(\V,\chi)=W_{\psi^{\kappa}}(\U,\mu).
\]
Furthermore, for every $\varphi\in \Sc(\V)$ and $\varphi_{0}$ defined
by equation \eqref{varphi0} we have
\[
W_{\varphi}=W_{\varphi_{0}}
\]
\end{prop}

\begin{pf} 
We define the map $j : \Sc(\V) \rightarrow W_{\psi^{\kappa}}(\V,\chi)$
by $j(\varphi)=W_{\varphi}$. The map $j$ is $\SLT(\F_{v})$ equivariant and it
is easy to see that it factors through the space  $\Sc(\V)_{\HH}$ of 
$\HH$ coinvariant. 

Since we know from Lemma \ref{par2} that $W_{\varphi}=W_{\varphi^{+}}$,
 the map $j$ actually factors through the $\GO(\V)$
coinvariants $\Sc(\V)_{\GO(\V)}$.

When $\chi=1$ the formula for 
$W_{\varphi}$ is
$$
W_{\varphi}(g) = \int_{T\backslash H} \omega(g)\varphi(h^{-1}\cdot x_{0})dh 
$$
If we let $\varphi$ be the characteristic function of some
neighborhood containing elements of norm $\kappa$ then it is easy to
see that $W_{\varphi}\neq 0$. Therefore $j\neq 0$.

On the other hand, 
let $ I(\eta)=\textup{Ind}_{B}^{\SLT}(\eta)$
where $\eta$ is the character \[
\begin{pmatrix} a & b \\
                                                0 & a^{-1}
                               \end{pmatrix} \mapsto
                               |a|^{\frac{3}{2}}\gamma(a,\psi^{-1}).
\]
There is an $\SLT(\F_{v})$ equivariant map $T:\Sc(\V)\rightarrow I(\eta)$
given by $T(\varphi)=(g \rightarrow
\omega(g)\varphi(0))$ which also factors through $\Sc(\V)_{\GO(\V)}$.
Let $R(\V)$ be the image of $T$. 
 A result of Rallis (\cite{R2} thm II.1.1) 
says that the map $T$ induces
an isomorphism $\Sc(\V)_{\GO(\V)} \simeq R(\V)$. It is also known
(e.g. Gelbart  \cite{G2}) that 
since $\V$ is a split space, then in fact $R(\V)\simeq I(\eta)$.
Furthermore, the composition series of $I(\eta)$ has length two, It has a unique
irreducible submodule and a unique irreducible quotient.

In
\cite{W2}, Waldspurger analyzed the two factors in the composition
series. His result is
that the unique quotient is $\Theta_{\U}({\bf 1})$: the local theta
correspondence of the trivial character from $\U$ (which Waldspurger
denotes by $r_{\psi,q}^{+}$.) This unique quotient has a
$\psi^{\kappa}$ Whittaker model whereas
the unique irreducible
submodule of $I(\eta)$ (which Waldspurger denotes by
$\tilde{\sigma}(|\,\cdot\,|^{\frac{1}{2}})$
has Whittaker models for all square classes of characters
$\psi^{\alpha}$ except
for $\psi^{\kappa}$ ( Lemma 3 on p 227 in \cite{W2}).

Since $j$ factors through $\Sc(\V)_{\GO(\V)} \simeq I(\eta)$
there is a map $\bar{j}: I(\eta) \rightarrow
W_{\psi^{\kappa}}(\V,\chi)$ and by considering the Whittaker models we
see that $\ker(\bar{j})=\tilde{\sigma}(|\,\cdot\,|^{\frac{1}{2}})$ and 
so $j(\Sc(\V))=W_{\psi^{\kappa}}(\V,\chi)$.
I.e. $j(\Sc(\V))$ is the unique irreducible quotient of $I(\eta)$
which is isomorphic to $\Theta_{\U}({\bf 1})$,
that is $j(\Sc(\V))=W_{\psi^{\kappa}}(\U,\mu)$.
\end{pf}

\section{Archimedean Local matching}
The argument splits into three cases handled separately in the next 
subsections.

\subsection{$\chi\circ\nu$ is trivial on $\GO(\V)(\R)$.}
Throughout this subsection, $\F_{v}$ will be isomorphic to $\R$ and we
will pick a particular realizations of the matrix algebra and division algebra
such that $\kappa=1$ when $\chi$ is nontrivial and $\kappa=-1$ when
$\chi$ is trivial. We will not mention these conditions in the
statements of the lemmas, but they should be understood to be stated
with
respect to these choices of $\kappa$. 
Since the equality we will establish is defined 
independently of a particular realization there is no loss of generality.

The cases included in this subsection are the case of $\chi =1$
( which by lemma \ref{ramtriv}, can only
exist when $\B$ is the matrix algebra) and the case $\chi \neq 1$ and 
$\B$ is the division algebra ( in which case  $\nu(\HH) \subset
\R^{+}$ so $\chi\circ \nu=1$). 
Note that in both cases, once we know that we can match all the
forms coming from the representation of $\SO(\V)(\R)$, by Lemma \ref{par2}, the
parity of all $\varphi\in\Sc(\V(\R))$ can be assumed to be even. Hence
only the trivial extension of 
$\chi\circ \nu$ to 
$\GO(\V)(\R)$ will have a nonzero lift. This is the reason for
describing both cases in this subsection as the case of $\chi\circ\nu$
trivial on $\GO(\V)(\R)$.

The matching in this case will be established by the matching of 
the Whittaker functions corresponding to the Gaussians and using 
the fact that the the space $\Theta(\HH,\chi)$ is generated by 
the image of the Gaussians. The matching of the Gaussians will 
depend on calculating their $K$-types and using the Iwasawa
decomposition to conclude that the two Whittaker functions agree on
all of $\SLT$.
The main argument is contained in the following theorem.

\begin{thm}\label{archimatch} 
Let $\chi$ be a the trivial character when $\B$ is a matrix algebra or
the nontrivial character when $\B$ is the division algebra.
Then $W_{\psi^{\kappa}}(\V,\chi)=W_{\psi^{\kappa}}(\U,\mu)$
and for any $\varphi\in \Sc(V)$ we have
and $W_{\varphi_{0}}=W_{\varphi}$. 
\end{thm}
\begin{pf}
As explained above, in this case the characters extend to the trivial
character on $\GO(\V)$. Therefore,  As mentioned in  \cite{KR}, it follows
from a result of Howe \cite{H2} that  $W_{\psi^{\kappa}}(\V,\chi)$ is generated
by the image of the Gaussian : $W_{\varphi}$.
The calculation that we carry out in this proof will show that if 
$\varphi$ is the Gaussian then 
$W_{\varphi_{0}}=W_{\varphi}$ and therefore 
$W_{\psi^{\kappa}}(\V,\chi)\subset W_{\psi^{\kappa}}(\U,\mu)$.
Therefore,
for any $\varphi$, there is 
some $f \in \Sc(\U)$ such that $W_{f}=W_{\varphi}$ 
and by Lemma \ref{bmatch} we know that $f$ must be $\varphi_{0}$.
This would also show that the formula for $\varphi_{0}$ defines a
Schwartz function in the restricted Schwartz space.

The second result is a consequence of Lemma
\ref{par2}, which implies
that the space $W_{\psi^{\kappa}}(\U,\mu)$
is the image of the irreducible
representation on even functions $\Sc(\U)^{+}$ so it
is also irreducible and $W_{\psi^{\kappa}}(\V,\chi)=W_{\psi^{\kappa}}(\U,\mu)$.

By the above argument, the calculation that remains is
to show that when $\varphi$ is the Gaussian we have
$W_{\varphi_{0}}=W_{\varphi}$. By Lemma \ref{weights}
we know that for all elements of
\[
K=\left\{ k_{\alpha} = \begin{pmatrix} \cos(\alpha) & -\sin(\alpha)\\
\sin(\alpha) & \cos(\alpha)\end{pmatrix} \right\}
\]
we have $\omega( k_{\alpha})\varphi=\exp(2\pi i m \alpha)\varphi$ and
$\omega( k_{\alpha})\varphi_{0}=\exp(2\pi i m \alpha)\varphi_{0}$ for
some $m$.

Using the Iwasawa decomposition
$\SLT=\widetilde{B}K$ we write an arbitrary element of $\SLT$ as $bk_{\alpha}$ for
some $\alpha$ and then 
\begin{align*}
W_{\varphi_{0}}(bk_{\alpha})
&=W_{\omega_{0}(k_{\alpha})\varphi_{0}}(b)\\
&=\exp(2\pi i m \alpha)W_{\varphi_{0}}(b)\\
&=\exp(2\pi i m \alpha)W_{\varphi}(b)\\
&=W_{\omega(k_{\alpha})\varphi}(b)\\
&=W_{\varphi}(bk_{\alpha})
\end{align*}
Where we have used the fact that by construction 
$W_{\varphi}(b)=W_{\varphi_{0}}(b)$.
\end{pf}

\begin{rem}
As we mentioned above we picked a realization of the algebras with 
$\nu(x_{0})=\kappa= 1,-1$ when $\chi$ is nontrivial or trivial,
respectively. By Lemma \ref{setm}, this means that we must choose the quadratic 
form on $\U$ to be
$x\rightarrow -x^{2}$ and $x\rightarrow x^{2}$ respectively .
\end{rem}

\begin{lem}\label{weights}
Let $\chi$ be trivial on $\nu(\HH)$ and let 
$\varphi \in \Sc(\V)$ be the Gaussian, 
then for some $m$
 $\omega( k_{\alpha})\varphi=\exp(2\pi i m \alpha)\varphi$ and
$\omega_{0}( k_{\alpha})\varphi_{0}=\exp(2\pi i m \alpha)\varphi_{0}$
\end{lem}
\begin{pf}
The Lie group $K$ is connected so its representation is determined by the action
of the Lie algebra. This Lie algebra is generated by 
the element $\zeta=\begin{pmatrix} 0 & -1 \\ 1 & 0 \end{pmatrix}$ so
it will be enough to calculate the actions of $\zeta$ to show that the 
actions of the whole algebra are the same.

We denote the Lie algebra action of $\omega$ and $\omega_{0}$ by
$\bar{\omega}$ and $\bar{\omega_{0}}$, respectively.
Instead of calculating the action of 
$\zeta$ directly, we write
$\zeta = e_{-}-e_{+}$ with
\[
e_{-}=\begin{pmatrix} 0 & 0 \\ 1 & 0 \end{pmatrix} \textup{ and }
e_{+}=\begin{pmatrix} 0 & 1 \\ 0 & 0 \end{pmatrix}
\]
and calculate the following actions
\[
\bar{\omega_{0}}(\zeta)=
\left. \frac{d}{dt}\right|_{t=0}\omega_{0}(\exp(te_{-})) + 
\left. \frac{d}{dt}\right|_{t=0}\omega_{0}(\exp(-te_{+}))
\]
and
\[
\bar{\omega}(\zeta)=
\left. \frac{d}{dt}\right|_{t=0}\omega(\exp(te_{-})) + 
\left. \frac{d}{dt}\right|_{t=0}\omega(\exp(-te_{+}))
\]
on $\varphi_{0}$ and $\varphi$ respectively. 
It is not hard to see that 
\[
\exp(-te_{+})=\begin{pmatrix} 1 & -t \\ 0 & 1 \end{pmatrix}=n(-t)
\]  and 
\[
\exp(te_{-})=\begin{pmatrix} 1 & 0 \\ t & 1 \end{pmatrix} =  w\exp(-te_{+})w^{-1}
\]
The rest of the proof is a lengthy but standard calculus
computation. Just as an illustration of the type of calculation involved we
 calculate the action on $\varphi$ when $\chi$ is trivial.
In this case the Gaussian is the function
\[
\varphi\left( \begin{pmatrix} u & v \\ w &- u \end{pmatrix}\right)=
\exp(- 2\pi|\epsilon|(u^{2}+\frac{v^{2}}{2}+ \frac{w^{2}}{2}))
\] 
and we calculate the first term in the action of $\zeta$.
\begin{align}
\bar{\omega}(-e_{+})\varphi(x)&=\left.
\frac{d}{dt}\right|_{t=0}\omega(\exp(-te_{+}))
\varphi(x)\nonumber \\
&=\left. \frac{d}{dt}\right|_{t=0}\omega(n(-t))\varphi(x)\nonumber \\
&=\left. \frac{d}{dt}\right|_{t=0}\exp(-t2\pi \epsilon i\nu(x))\varphi(x) 
\nonumber \\
&= -2\pi \epsilon i\nu(x)\varphi(x) \label{e+1}
\end{align}
We now calculate the second term in the action of $\zeta$.
\begin{align*}
\bar{\omega}(e_{-})\varphi(x)&=
\left. \frac{d}{dt}\right|_{t=0}\omega(\exp(te_{-}))\varphi(x) \\
&=\left. \frac{d}{dt}\right|_{t=0}\omega(wn(-t)w^{-1})\varphi(x)\\
&=\left. \frac{d}{dt}\right|_{t=0}\gamma(\psi\circ\V)\omega(wn(-t)\varphi(x)
\end{align*}
The last step uses a Lemma (which is also omitted) which says that 
the Gaussian is an eigenfunction of $\omega(w)$ with eigenvalue 
$\gamma(\psi\circ \V)^{-1}$. 
\begin{align*}
\bar{\omega}(e_{-})
\varphi(x)
&=\left. \frac{d}{dt}\right|_{t=0}\int_{\V}\omega(n(-t))\varphi(y)\psi(-(x,y))dy
\end{align*}
To continue the calculation we let 
\[ x = \begin{pmatrix} \alpha & \beta \\ \gamma &- \alpha \end{pmatrix}
\textup{ and }
y = \begin{pmatrix} u & v \\ w &-u \end{pmatrix}
\]
and then
\begin{align*}
\bar{\omega}(e_{-})\varphi(x)
=\left. \frac{d}{dt}\right|_{t=0}
\int\int\int&\exp(2\pi \epsilon i t(u^{2} + vw)) \\
&\exp(-2\pi |\epsilon|(u^{2}+\frac{v^{2}}{2}+ \frac{w^{2}}{2}))\\
&\exp(2\pi \epsilon i (2u\alpha + \gamma v + \beta w))du\,dv\,dw 
\end{align*}
Using standard formulas from calculus and integrating with respect to \
the self dual measure
\begin{align*}
\bar{\omega}(e_{-})\varphi(x)=&\frac{i\textup{sgn}(\epsilon)
\exp(-2\pi|\epsilon|(\alpha^{2}+\frac{\beta^{2}}{2}+\frac{\gamma^{2}}{2}))}
{2}\\ 
&+ 2\pi\epsilon i(-\alpha^{2}-\beta\gamma)
\exp(-2\pi|\epsilon|(\alpha^{2}+\frac{\beta^{2}}{2}+\frac{\gamma^{2}}{2}))
\end{align*}
so that
\begin{equation}\label{e-1}
\bar{\omega}(e_{-})\varphi(x)
=\frac{i\textup{sgn}(\epsilon)}{2}\varphi(x) + 2\pi i \epsilon \nu(x)\varphi(x)
\end{equation}
Adding up equation \eqref{e-1} and \eqref{e+1} we get the action of
$\zeta$ on the Gaussian
\begin{align}\label{zeta-gaussian}
\bar{\omega}(\zeta)\varphi(x)&=
\bar{\omega}(e_{-})\varphi(x) + \bar{\omega}(-e_{+})\varphi(x)\\
&=\frac{i\textup{sgn}(\epsilon)}{2}\varphi(x)
\end{align}
\end{pf}

\subsection{$\chi\circ \nu$ is nontrivial on $\HH(\R)$.}

The remaining archimedean case is the case of a nontrivial $\chi$ and 
$\B$ is a
matrix algebra. The argument that will show matching is similar to
the one for the other real cases. We will use the fact that the
representation is known to be irreducible and the matching of one
function
to conclude that the Whittaker spaces of $W_{\psi^{\kappa}}(\V,\chi)$ 
and $W_{\psi^{\kappa}}(\U,\mu)$
coincide.

Throughout this section $\F_{v}\simeq \R$, $\chi$ is nontrivial and $\B\simeq\M_{2}$.

\begin{thm} \label{archimatch2}
Let $\chi$ be the nontrivial character on $\R$ then 
\[W_{\psi^{\kappa}}(\V,\chi)=W_{\psi^{\kappa}}(\U,\mu)\]
and for any 
$\varphi\in \Sc(\V)$ $W_{\varphi}=W_{\varphi_{0}}$.
\end{thm}
\begin{pf}
Since in this case the quadratic space $\V$ is split with signature (2,1)
we can use a theorem of  Zhu, (1.2 in \cite{Z}) which says that in
this case, when $\chi=\textup{sgn}$,
$W_{\psi^{\kappa}}(\V,\textup{sgn})=W_{\psi^{\kappa}}(\V,{\bf 1})$.
By Lemma \ref{archimatch} we
know that $W_{\psi^{\kappa}}(\V,{\bf 1})$ 
is irreducible and therefore so is $W_{\psi^{\kappa}}(\V,\chi)$.

Since $W_{\psi^{\kappa}}(\V,\chi)$ is irreducible and since we also know that 
$W_{\psi^{\kappa}}(\U,\mu)$ is irreducible it will be enough to show that
their Whittaker spaces intersect to show that they coincide.

We will show that when $\varphi$ is the gaussian then
$W_{\varphi}=W_{\varphi_{0}}$, this would mean that the
Whittaker spaces are equal and by Lemma \ref{bmatch} it follows that
for any $\varphi \in \Sc(\V)$ we have, $W_{\varphi}=W_{\varphi_{0}}$.
It will also follow that $\varphi_{0}$ is actually a function in the
restricted Schwartz space of $\Sc(\R)$.
As in the proof of theorem \ref{archimatch} it will be enough to show
that when $\varphi$ is the gaussian then $\varphi$ and $\varphi_{0}$
are of the same K-type. This is done in Lemma \ref{weights2}
and completes the argument.
\end{pf}

\begin{lem}\label{weights2}
Let $\chi$ be a nontrivial character, let $\B$ be the matrix algebra
and let $\varphi \in \Sc(\V)$ be
the gaussian character, then for some $m$, 
$\omega( k_{\alpha})\varphi=\exp(2\pi i m \alpha)\varphi$ and
$\omega( k_{\alpha})\varphi_{0}=\exp(2\pi i m \alpha)\varphi_{0}$.
\end{lem}
\begin{pf}
omitted\end{pf}

\subsection{The case of $\GO(\V)(\C)$}
The local data in the complex case is as simple as possible. Since
every complex number is a square, the quadratic character
$\chi=\chi_{v}$ is trivial. For the same reason, the quaternion
algebra is ramified i.e. $\textup{Inv}(\B_{v})=1$.

The correspondence in this case can be deduced from known facts
without any additional calculation. In \cite{AB} (Proposition 2.1)
Adams and Barbasch
show that the theta lift of the trivial character from $\GO(\V)(\C)$
to $\SL_{2}(\C)$ is isomorphic to the irreducible representation 
$\omega_{+}$ which is the Weil representation of $\SL_{2}(\C)$ on
$\Sc(\U)(\U)^{+}$- the space of even Schwartz functions.

If $I'_{\V}(\chi\otimes{\bf 1},\varphi)$ is the automorphic form in
the theta lift of the
character $\chi\circ\nu$ extended trivially from $\SO(\V)(\C)$ to
$\GO(\V)(\C)$ defined similarly to $I_{\V}(\chi,\varphi)$ then
 $I'_{\V}(\chi\otimes{\bf 1},\varphi)=2I_{\V}(\chi,\varphi)$. 
Therefore, $\Theta_{V}(\chi)\simeq \omega^{+}\simeq \Theta_{\U}(\mu)$ 
and to conclude the
matching of functions we need to show that the spaces of functions
$\Theta_{\V}(\chi)$ 
and
$\Theta_{\U}(\mu)$ actually coincide. This follows from the
multiplicity one of spherical representations of $\SL_{2}(\C)$.
The fact that the representations are spherical is verified by the
existence of the gaussian functions in both representations.
Therefore, the matching for the actual function is given by the
formula for $\varphi_{0}$ \eqref{varphi0}. 

\begin{thebibliography}{}
\expandafter\ifx\csname url\endcsname\relax
  \def\url#1{\texttt{#1}}\fi
\expandafter\ifx\csname urlprefix\endcsname\relax\def\urlprefix{URL }\fi

\end{thebibliography}


\begin{thebibliography}{00}




\bibitem{AB}
J~Adams and D~Barbasch, \emph{Reductive dual pair correspondence for complex
  groups}, Jour Functional Analysis. \textbf{132} (1995), 1--42.

\bibitem{B}
D.~Bump, \emph{Automorphic forms and representations}, Cambridge University
  Press, Cambridge, 1997. 

\bibitem{Z}
Zhu Chen-Bo, \emph{Representations with scalar k-types and applications},
  Israel Jour. Math. \textbf{135} (2003), no.~1, 111--124.

\bibitem{GPS}
S.~Gelbart and I.~I. Piatetski-Shapiro, \emph{Distinguished representations and
  modular forms of half-integral weight}, Invent. Math. \textbf{59} (1980),
  no.~2, 145--188. 

\bibitem{G2}
Stephen~S. Gelbart, \emph{Weil's representation and the spectrum of the
  metaplectic group}, Springer-Verlag, Berlin, 1976, Lecture Notes in
  Mathematics, Vol. 530. 

\bibitem{H2}
R.~Howe, \emph{Transcending classical invariant theory}, J. AMS. \textbf{2}
  (1989), no.~3, 535--552.

\bibitem{KU}
S.~Kudla, \emph{On the theta correspondence}, unpublished notes, College Park,
  1996.

\bibitem{KR}
Stephen~S. Kudla and Stephen Rallis, \emph{Degenerate principal series and
  invariant distributions}, Israel J. Math. \textbf{69} (1990), no.~1, 25--45.
  

\bibitem{SP2}
Schulze-Pillot Rainer, \emph{Thetareihen positiv definiter quadratischer
  formen.}, Invent. Math \textbf{75} (1984), 283--299.

\bibitem{R2}
S.~Rallis, \emph{On the {H}owe duality conjecture}, Compositio Math.
  \textbf{51} (1984), no.~3, 333--399. 

\bibitem{RA}
R.~Ranga~Rao, \emph{On some explicit formulas in the theory of {W}eil
  representation}, Pacific J. Math. \textbf{157} (1993), no.~2, 335--371.
 

\bibitem{SE}
J.-P. Serre, \emph{A course in arithmetic}, Springer-Verlag, New York, 1973,
  Translated from the French, Graduate Texts in Mathematics, No. 7.

\bibitem{VI}
Marie-France Vign{\'e}ras, \emph{Arithm\'etique des alg\`ebres de quaternions},
  Lecture Notes in Mathematics, vol. 800, Springer, Berlin, 1980.
 

\bibitem{W1}
J.-L. Waldspurger, \emph{Correspondance de {S}himura}, J. Math. Pures Appl. (9)
  \textbf{59} (1980), no.~1, 1--132.

\bibitem{W2}
J.-L. Waldspurger, \emph{Correspondances de {S}himura et quaternions}, Forum Math.
  \textbf{3} (1991), no.~3, 219--307. 


\end{thebibliography}
\bibliographystyle{elsart-num}




Kobi Snitz,         
                  The Centre for Advanced Studies in Mathematics 
                  Ben Gurion University of the Negev 
                  P.O.B. 653 
                  Be'er Sheva 84105 
                  ISRAEL  
\end{document}